\documentclass[4p,12pt]{article}
\usepackage{amssymb}
\usepackage{amsmath,amsthm}
\usepackage{amsfonts}
\usepackage{graphicx}
\usepackage{graphics}

\newcommand{\be}{\begin{equation}}
\newcommand{\ee}{\end{equation}}
\newcommand{\bes}{\begin{equation*}}
\newcommand{\ees}{\end{equation*}}

\begin{document}

\date{\small\textsl{\today}}
\title{The dual reciprocity boundary elements method for one--dimensional nonlinear parabolic partial differential equations}
\author{ \large Peyman Alipour \footnote{Corresponding author.\newline {\em  E-mail
addresses:} palipour@stevens.edu. \newline} $\vspace{.2cm} $ \large \\
\small{\em School of Business, Stevens Institute of Technology, Hoboken, NJ, USA}\vspace{-1mm}} \maketitle

  \vspace{-0.1cm}

{\bf Abstract}:
This article describes a numerical method based on the dual reciprocity boundary elements method (DRBEM)
for solving some well--known nonlinear parabolic partial differential equations (PDEs).
The equations include the classic and generalized Fisher's equations, Allen--Cahn equation, Newell--Whithead equation,
FitzHugh--Nagumo equation and generalized
FitzHugh--Nagumo equation with time--dependent coefficients.
The concept of the dual reciprocity is used to convert the domain integral to the boundary that leads to an integration free method.
We employ the time stepping scheme to approximate the time derivative, and the linear radial basis functions (RBFs) are used as
approximate functions in presented method. The nonlinear terms are treated iteratively within each time step.
The developed formulation is verified in some numerical test examples. The results of numerical experiments are compared with analytical solution to confirm the accuracy and efficiency of the presented scheme. \\



\section{Introduction}
The dual reciprocity boundary element method (DRBEM) is a modern numerical scheme, which has enjoyed increasing popularity and it has become one of the most general and effective numerical method for solving different engineering problems~\cite{ang,DRMbook,Pozrikidis}. Generally, the DRBEM is known as a numerical tool for solving two or higher dimensional partial differential equations (PDEs). However, some authors employed this approach for the numerical solution of some one--dimensional PDEs. For examples, hyperbolic telegraph equation \cite{ghesmati_OneDim_telegraph_DRM},
nonlinear Klein–-Gordon equation \cite{ghesmati_OneDim_KG_DRM}, nonlinear sine--Gordon equation \cite{dav_OneDim_SG_DRM},
Cahn–-Hilliard equation \cite{dav_OneDim_CH_DRM} and advection--diffusion equation \cite{shiva}
have been solved with the one--dimensional DRBEM approach.\\
In this study, a numerical solution based on the DRBEM is applied for solving some well--known one--dimensional nonlinear parabolic PDEs. The idea behind this approach comes from the classic DRBEM introduced by Brebbia and Nardini
\cite{firstDRM} and Partridge and Brebbia \cite{DRMbook} for solving higher order dimensional problems and to expand the inhomogeneous
and nonlinear terms in terms of their values at the nodes which lie in the domain of problem.
The inhomogeneous term is approximated by interpolation in terms of some well--known functions $\phi(r)$,
called radial basis functions (RBFs), where $r$ is the distance between a source point and the field point.\\

Our concern in the current work is to present a numerical method based on the DRBEM for solving the following nonlinear parabolic PDEs
\begin{equation}\label{model}
u_t+\nu(t)u_x-\mu(t)u_{xx}-\eta(t)F(u)=0, \ \ \ \ (x,t) \in [a,b] \times [0,T],
\end{equation}
subject to the initial condition
\begin{equation}
u(x,0)=f(x), \ \ \ x \in [a,b],
\end{equation}
and the boundary conditions
\begin{equation}
u(a,t)=g_1(t), \ \ \ \ \ \ \ \ u(b,t)=g_2(t), \ \  t \in [0,T],
\end{equation}
where $\nu(t), \mu(t)$ and $\eta(t)$ are arbitrary real--valued functions of $t$. $F(u)$ in Eq. (\ref{model}), written as $F_l(u)+F_n(u)$, where $F_l$ and $F_n$ denote the linear and nonlinear parts of $F$, respectively. Eq. (1) for different values of $\nu, \mu, \eta$ and $F$ yields the following well--known problems:
\subsection*{\small{Case 1: The Fitzhugh--Nagumo and real Newell--Whitehead equations}}
If we set $\nu(t)=\mu(t)=\eta(t)=1$ and $F(u)=u(1-u)(\rho-u)$ where $0 \leq \rho \leq 1$, Eq. (1) deduces to classic Fitzhugh--Nagumo equation. In addition, if in Fitzhugh--Nagumo equation $\rho$ takes the value $-1$ then the classic Fitzhugh--Nagumo equation deduce to real Newell--Whitehead equation. The Fitzhugh--Nagumo equation has been derived by Fitzhugh \cite{F} and Nagumo et al. \cite{N}. Also population genetics \cite{populationgenetic} is another area of application of the Fitzhugh--Nagumo equation.
Kawahara and Tanaka~\cite{Kawahara}, Nucci and Clarkson~\cite{Nucci}, Li and Guo~\cite{Li}, Abbasbandy \cite{abbasbandy} have found some new solution of the Fitzhugh--Nagumo equation using, respectively, Hirota method, Jacobi elliptic function, first integral method and homotopy analysis method.
The Haar wavelet method \cite{Hariharan}, Pseudospectral methods \cite{fakhar}, the homotopy analysis method \cite{homotopy} are some numerical approach that have been applied to present the approximate solution of the Fitzhugh--Nagumo equation. Meanwhile the authors of \cite{Shih} proposed the approximate conditional symmetry method to determine approximate solutions of Fitzhugh--Nagumo equation.
\subsection*{\small{Case 2: The generalized Fitzhugh--Nagumo equation}}
Let $\nu(t), \mu(t)$ and $\eta(t)$ be an arbitrary function of $t$ and let $F(u)=u(1-u)(\rho-u)$. Then Eq. (1) deduces to generalized Fitzhugh--Nagumo with
time--dependent coefficients and linear dispersion term equation.
The authors of \cite{wazwazFN} derived new variety of soliton solutions using specific solitary wave ansatz and the tanh method for this equation. Meanwhile Bhrawy \cite{Bhrawy} employed Jacobi-–Gauss-–Lobatto collocation method for the numerical solution of this equation.
\subsection*{\small{Case 3: The Fisher's equation}}
The Fisher's equation achieved when $\nu(t)=\mu(t)=\eta(t)=1$ and $F(u)=u(1-u)$. As mentioned in \cite{wazwazFisher}, Fisher proposed
such equation as a model for the propagation of a mutant gene, with $u$ denoting the density of an advantageous. This equation
is encountered in chemical kinetics and population dynamics which include problems such as nonlinear evolution of a
population in a one--dimensional habitat, neutron population in a nuclear reaction. Moreover, the same equation occurs
in logistic population growth models, flame propagation, neurophysiology, autocatalytic chemical reactions,
and branching Brownian motion processes. Sinc collocation method \cite{Khaled}, B--spline Galerkin method \cite{dag,Mittal2}, wavelet Galerkin method \cite{Mittal}, moving mesh method \cite{Qiu}, finite element methods \cite{Ro,Tang} and finite difference \cite{Mickens} are some numerical techniques that have been applied for the numerical solution of Fisher's equation. Meanwhile Wazwaz \cite{wazwazFisher,wazwaz} have been found some exact solution for this equation with tanh-coth and Adomian decomposition methods.
\subsection*{\small{Case 4: The Allen-Cahn and generalized Fisher's equations}}
If we set $\nu(t)=\mu(t)=\eta(t)=1$ and $F(u)=u(1-u^{\alpha})$ in Eq. (\ref{model}), the generalized Fisher's equation will be obtained. In addition, Eq. (1) for $\alpha=$2 and $\alpha>2$ are called the Allen-Cahn equation and the generalized Fisher's equation. Some useful numerical solutions of this equation are \cite{Fen,25,chino,HariharanCH}. In addition, \cite{wazwaz} provides some exact solutions of these equations.
\subsection{The organization of the current paper}
This article is organized as follows: In Section 2, the discretized version of equation is obtained and an iterative scheme based on finite difference scheme is described for the time derivative. In Section 3, numerical results of some nonlinear parabolic PDEs have been presented and the obtained results are compared with the exact solutions. Section 4 ends this report with a brief conclusion.
\section{The proposed method}
Suppose $ G_i=G(x,x_i) $ is the fundamental solution of the one--dimensional Laplace operator based on the source point $x_i$,
i.e. \\
\begin{equation}\label{deltadirak}
\frac{\partial^2 G}{\partial x^2} (x,x_i)=\delta(x,x_i),
\end{equation}
where $x$ is field point and  $\delta$ is Dirac delta function. The fundamental solution and its derivative are
given as follows
\cite{fundamentalsolution}:
\begin{equation}
\begin{array}{l}
G(x,x_i)=\frac{1}{2}|x-x_i|, \\
\\
G_{,x}(x,x_i)=G'(x,x_i)=\frac{1}{2}sgn(x-x_i),
\end{array}
\end{equation}
where the symbol $sgn$ denotes the signum function.\\
Consider Eq. (\ref{model}) as follows:
\begin{equation}
u_{xx}=\frac{1}{\mu(t)} u_t+\frac{\nu(t)}{\mu(t)} u_x-\frac{\eta(t)}{\mu(t)} \{F(u)\}:=b(x,t).
\end{equation}
Multiplying the above equation by $G_i$, taking integration over $[a,b]$ and applying the integration by parts,
we get the following integral form:
\begin{equation}\label{intb1f}
[G_iu_x]_a^b-[G_i'u]_a^b+c_iu_i=\int_a^b b(x,t)G_idx,
\end{equation}
where similarly to 2D formulation, $c_i$ takes the values $1$ and $1/2$ when the source point is located in
domain and on the boundary (the points $a$ and $b$), respectively \cite{DRMbook,adeliz}. \\
The  domain integral on the right hand side of Eq. (\ref{intb1f}) still remains in the boundary elements formulation.
This integral can be evaluated by dividing the domain into cells \cite{rbfdrm}.
The motivation behind DRBEM is to avoid this procedure by transforming the domain integral to an equivalent boundary integral equation.
This can be achieved by approximating the function $ b(x,t)$ in terms of radial basis functions (RBFs)
at some chosen number in
$[a,b]$ as
\begin{equation}
a=x_1 <x_2 <\cdots<x_{N-1}<b=x_N.
\end{equation}
So the function  $ b $  can be  expressed  as
\begin{equation} \label{bf}
b= \sum_{j=1} ^{N} \alpha_j(t) \phi_j,
\end{equation}
where $\alpha_j(t), j=1,...,N$, are the corresponding interpolating
coefficients, $\phi_j$  represents the interpolation function, $\phi$, from a field node to source node, i.e.
\begin{equation*}
\phi_j:=\phi(|x-x_j|), \ \ x, x_j \in \mathcal{D},
\end{equation*}
where $|x-x_j|$ denotes the distance between $x$ and $x_j$.
The essential feature in DRBEM is to express $\phi_j$,  which is a function of $r_j$,
as a Laplacian of another function $\psi _j$. Thus $\psi _j$ is chosen as the solution to \cite{rbfdrm}
\begin{equation}\label{nabla2rbf}
\frac{\partial^2 \psi _j}{\partial x^2} =\phi_j.
\end{equation}
In this paper, we will use linear RBFs as
\begin{equation}
\phi_j=1+r_j.
\end{equation}
The function $\psi _j$ is easily determined as
\begin{equation}
\psi _j=\frac{1}{2} r_j^2+\frac{1}{6} r_j^3.
\end{equation}
With substitution expansion (\ref{bf}) for $b(x,t)$, applying the integration by part once again,
the domain integral in the right hand side of Eq. (\ref{intb1f}) reduces to a boundary integral equation
\begin{equation}\label{rhsb}
\int_a^b b(x,t) G_i dx=\sum_{j=1}^{N} \{ [G_i\psi_j']_a^b-[G_i'\psi_j]_a^b+c_i\psi_{ij} \} \alpha _j(t),\\
\end{equation}
where  $\psi_{ij} $  is  the  value  of  the  function  $\psi_j$ at the $i$th source point.
So from Eqs. (\ref{intb1f}) and (\ref{rhsb}) the following boundary integral equation can be achieved
\begin{equation}\label{bie}
[G_iu_x]_a^b-[G_i'u]_a^b+c_iu_i=\sum_{j=1}^{N} \{ [G_i\psi_j']_a^b-[G_i'\psi_j]_a^b+c_i\psi_{ij} \} \alpha _j(t).
\end{equation}
Imposing all the source points to satisfy Eq. (\ref{bie}) yields the following matrix form:
\begin{equation}\label{313}
\begin{array}{l}
\bold{L}\left[
          \begin{array}{c}
            u_x(a) \\
            u_x(b) \\
          \end{array}
        \right]-\bold{H}\left[
                          \begin{array}{c}
                            u(a) \\
                            u(b) \\
                          \end{array}
                        \right]+\left[
                                  \begin{array}{c}
                                    \frac{1}{2}u(a) \\
                                    \bold{u}_{in} \\
                                    \frac{1}{2}u(b) \\
                                  \end{array}
                                \right],
\\
=\displaystyle \sum_{j=1}^{N} \{ \bold{L}\left[
          \begin{array}{c}
             \psi_j'(a) \\
             \psi_j'(b) \\
          \end{array}
        \right]-\bold{H}\left[
                          \begin{array}{c}
                          \psi_j(a) \\
                          \psi_j(b) \\
                          \end{array}
                        \right]+\left[
                                  \begin{array}{c}
                                    \frac{1}{2}\psi_j(a) \\
                                    \psi_j(x_1) \\
                                    \vdots\\
                                    \psi_j(x_{N-1})\\
                                    \frac{1}{2}\psi_j(b) \\
                                  \end{array}
                                \right] \} \alpha_j(t),\\
\end{array}
\end{equation}
where $\bold{u}_{in}=[u(x_2),\cdots,u(x_{N-1})]^T$, $u_x(a)=u_x(a,t), u_x(b)=u_x(b,t), u(x_j)=u(x_j,t), j=1,...,N$
and $\bold{L}$ and $\bold{H}$ take the following form
\begin{equation*}
\bold{L}=\left[
           \begin{array}{cc}
             -G_1(a) & G_1(b) \\
             -G_2(a) & G_2(b) \\
             \vdots & \vdots \\
             -G_N(a) & G_N(b) \\
           \end{array}
         \right], \ \ \ \ \ \ \ \bold{H}=\left[
           \begin{array}{cc}
             -G_1'(a) & G_a'(b) \\
             -G_2'(a) & G_1'(b) \\
             \vdots & \vdots \\
             -G_N'(a) & G_N'(b) \\
           \end{array}
         \right].
\end{equation*}
If each of the vectors
\begin{equation*}
\begin{array}{l}
\left[
          \begin{array}{c}
            \psi_j'(a) \\
            \psi_j'(b) \\
          \end{array}
        \right], \ \ \ \ \left[
                          \begin{array}{c}
            \psi_j(a) \\
            \psi_j(b) \\
                          \end{array}
                        \right] \ \ and \ \ \ \ \left[
                                  \begin{array}{c}
                                    \frac{1}{2}\psi_j(a) \\
                                    \psi_j(x_2) \\
                                    \vdots\\
                                    \psi_j(x_{N-1})\\
                                    \frac{1}{2}\psi_j(b) \\
                                  \end{array}
                                \right],
\end{array}
\end{equation*}
are considered to be one column of the matrices $ \Psi_x$, $ \Psi $ and $ \widetilde{\Psi}$, respectively, Eq. (\ref{313})
takes the following matrix form
\begin{equation}\label{31}
\begin{array}{l}
\bold{L}\left[
          \begin{array}{c}
            u_x(a) \\
            u_x(b) \\
          \end{array}
        \right]-\bold{H}\left[
                          \begin{array}{c}
                            u(a) \\
                            u(b) \\
                          \end{array}
                        \right]+\left[
                                  \begin{array}{c}
                                    \frac{1}{2}u(a) \\
                                    \bold{u}_{in} \\
                                    \frac{1}{2}u(b) \\
                                  \end{array}
                                \right] =[\bold{L} \Psi_x - \bold{H}\Psi+\widetilde{\Psi}]\alpha.
                                \end{array}
\end{equation}
On the other hand Eq. (\ref{bf}) can be written in the following matrix form
\begin{equation}
\bold{\Phi} \alpha=\bold{b},
\end{equation}
where $\Phi_{ij}$ represents the value of the function $\phi_j$ at source point $x_i$ by $\phi_{ij}$ for $i =1, 2,...,N,$
and vector $\bold{b}$ takes the following form
\begin{equation}
\bold{b}=\left[
  \begin{array}{c}
    b(x_1,t) \\
    \vdots \\
    b(x_N,t) \\
  \end{array}
\right],
\end{equation}
where
\begin{equation}
b(x_j,t)=\frac{1}{\mu(t)} u_t(x_j,t)+\frac{\nu(t)}{\mu(t)} u_x(x_j,t)-
\frac{\eta(t)}{\mu(t)}(u(x_j,t)\{ F(u(x_j,t))\}.
\end{equation}
Now Eq. (\ref{31}) constitutes a nonlinear system of $N$ equations in $N$ unknown functions of $t$.
This system is solved approximately using the iterative scheme based on the implicit finite difference technique as follow:
\begin{equation}\label{m2}
\begin{array}{l}
\bold{L}\left[
          \begin{array}{c}
            u^n_x(a) \\
            u^n_x(b) \\
          \end{array}
        \right]-\bold{H}\left[
                          \begin{array}{c}
                            u^n(a) \\
                            u^n(b) \\
                          \end{array}
                        \right]+\left[
                                  \begin{array}{c}
                                    \frac{1}{2}u^n(a) \\
                                    \bold{u}^n_{in} \\
                                    \frac{1}{2}u^n(b) \\
                                  \end{array}
                                \right]= \bold{D}\bold{\Phi}^{-1} \bold{b},\\
\end{array}
\end{equation}
where $\bold{D}=[\bold{L} \Psi_x - \bold{H}\Psi+\widetilde{\Psi}]$ and $\bold{b}$ takes the following form
\begin{equation}\label{bfinal}
\bold{b}=\frac{\bold{u}^n-\bold{u}^{n-1}}{\tau \mu(t_n)}+\frac{\nu(t_n)}{\mu(t_n)} \bold{u}^n_x- \frac{\eta(t_n)}{\mu(t_n)}\{ F_l(\bold{u}^n)+F_l(\widetilde{\bold{u}})\}.
\end{equation}
On the other hand the solution of Eq. (\ref{model}) can be approximated as follows \cite{TEZboundaryelement}
\begin{equation}\label{uap}
u\simeq \sum_{j=1}^N \beta_j \phi_j,
\end{equation}
or in matrix form
\begin{equation}\label{um}
\bold{u}=\bold{\Phi} \beta(t).
\end{equation}
Therefore from Eqs. (\ref{uap}) and (\ref{um}) we can write
\begin{equation} \label{ux}
\bold{u}_x=\bold{\Phi}_x \beta(t)=\bold{\Phi}_x \bold{\Phi}^{-1} \bold{u}.
\end{equation}
So by substituting (\ref{ux}) in (\ref{bfinal}) and substituting the result in (\ref{m2}) and separating the known quantities from the unknown quantities,
the value of $\bold{u}_{in}^n$ can be obtained by solving a nonlinear system of equation.
To avoid solving the nonlinear system of equations, the following iterative algorithm has been proposed \cite{shokri2,davLBIE,ahMLS,shokriPC}.\\
\subsection{The predictor-corrector scheme}
For dealing with the non--linearity, in time level $n$ at first put
\begin{equation}
\widetilde{\bold{u}}=\bold{u}^{n-1}.
\end{equation}
With this substitution, Eq. (\ref{m2}) is solved as a system of linear algebraic equations for unknown $\bold{u}^n=\bold{u}^{n,0}$.
Recompute
\begin{equation}
\widetilde{\bold{u}}=\bold{u}^{n,0}.
\end{equation}
Now Eq. (\ref{m2}) is solved using the new $\widetilde{\bold{u}}$ for unknown $\bold{u}^{n,l}$.
We are at time level $n$ yet, and iterate between calculating $\widetilde{\bold{u}}$ and computing the approximation values of
the unknown
$\bold{u}^{n,l}$ and putting
\begin{equation}
\widetilde{\bold{u}}=\bold{u}^{n,l},
\end{equation}
until the unknown quantity converges to within a prescribed number of the significant figures.
In this paper, we will use the following condition for stopping the iterations in each time level:
\begin{equation}
\parallel \bold{u}^{n,l}-\bold{u}^{n,l-1}  \parallel_{\infty} \leq \epsilon,
\end{equation}
where $\epsilon$ is a fixed number. When this condition is satisfied we put
\begin{equation}
\widetilde{\bold{u}}=\bold{u}^{n,l},
\end{equation}
and go on to the next time level. This process is iterated, until reaching to the desirable time $t$.
\section{Numerical simulations}

To access both important parameters in numerical solution of problems, the accuracy and the applicability of the procedure
described in the previous section, some test examples are considered. In the following test problems we will use the
$L_{\infty}$ and the root--means--squares ($RMS$) errors,
as defined bellow, to report the errors
\begin{equation}\label{errors}
\begin{array}{l}
L_{\infty}-error=\displaystyle{\max_{1 \leq j \leq N}|e_j|}, \\
RMS-error=\displaystyle {\sqrt{\frac{1}{N} \sum_{j=1}^N |e_j|^2}}, \\
\end{array}
\end{equation}
where
\begin{equation*}
\begin{array}{l}
e=\displaystyle {u_{exact}-u_{approximate}},\\
e_j=\displaystyle {(u_{exact}-u_{approximate})_j}, \ \ \ \ j=2,...,N-1. \\
\end{array}
\end{equation*}
Also in this part we assume $\epsilon=10^{-10}$.
\subsection{Example 1}

Consider the following Fitzhugh–-Nagumo equation
\begin{equation}\label{fn}
u_t=u_{xx}-u(1-u)(\rho-u), \ \ \ (x,t) \in [a,b] \times [0,T],
\end{equation}
subject to the initial condition
\begin{equation}
u(x,0)=\frac{1}{2}+\frac{1}{2} \tanh(\frac{x}{2 \sqrt{2}}), \ \ \ x \in [a,b],
\end{equation}
and the following boundary conditions
\begin{equation}
\begin{array}{l}
u(a,t)=\displaystyle(\frac{1}{2}+\frac{1}{2} \tanh(\frac{1}{2 \sqrt{2}}(a-\frac{2\rho-1}{\sqrt{2}}t))), \\
\\
u(b,t)=\displaystyle(\frac{1}{2}+\frac{1}{2} \tanh(\frac{1}{2 \sqrt{2}}(b-\frac{2\rho-1}{\sqrt{2}}t))). \\
\end{array}
\end{equation}
The exact solution of Eq. (\ref{fn}) is given by \cite{Bhrawy,wazwaz}
\begin{equation}
u(x,t)=\displaystyle(\frac{1}{2}+\frac{1}{2} \tanh(\frac{1}{2 \sqrt{2}}(x-\frac{2\rho-1}{\sqrt{2}}t))), \\
\end{equation}
The error norms defined by (\ref{errors}) on grids for three constant $\tau $ with an increasing number of nodes are presented in Table 1. The results reveal that the error decreases when the number of nodes increases.
\begin{table}[h]
\caption{The obtained estimate errors for Ex. 1}
\begin{center}
\begin{tabular}{cccc}
  \hline
  \hline
  $h$ & $\tau$ & $L_{\infty}$ & $RMS$ \\
  \hline
  \hline
  \footnotesize{$1/4$}  & \footnotesize{$1/500$}& \footnotesize{$2.8473\textbf{E}-05$} & \footnotesize{$1.9553\textbf{E}-05$} \\
  \footnotesize{$1/8$}  & \footnotesize{$1/500$}& \footnotesize{$8.0477\textbf{E}-06$} & \footnotesize{$5.2319\textbf{E}-06$} \\
  \footnotesize{$1/16$} & \footnotesize{$1/500$}& \footnotesize{$2.7413\textbf{E}-06$} & \footnotesize{$1.7839\textbf{E}-06$} \\
  \footnotesize{$1/32$} & \footnotesize{$1/500$}& \footnotesize{$1.4138\textbf{E}-06$} & \footnotesize{$9.3950\textbf{E}-07$} \\
  \footnotesize{$1/64$} & \footnotesize{$1/500$}& \footnotesize{$1.0823\textbf{E}-06$} & \footnotesize{$7.3079\textbf{E}-07$} \\
  \hline
  \hline
  \footnotesize{$1/4$}  & \footnotesize{$1/1000$}& \footnotesize{$2.7938\textbf{E}-05$} & \footnotesize{$1.9159\textbf{E}-05$} \\
  \footnotesize{$1/8$}  & \footnotesize{$1/1000$}& \footnotesize{$7.5096\textbf{E}-06$} & \footnotesize{$4.8528\textbf{E}-06$} \\
  \footnotesize{$1/16$} & \footnotesize{$1/1000$}& \footnotesize{$2.2037\textbf{E}-06$} & \footnotesize{$1.4083\textbf{E}-06$} \\
  \footnotesize{$1/32$} & \footnotesize{$1/1000$}& \footnotesize{$8.7637\textbf{E}-07$} & \footnotesize{$5.6270\textbf{E}-07$} \\
  \footnotesize{$1/64$} & \footnotesize{$1/1000$}& \footnotesize{$5.4445\textbf{E}-07$} & \footnotesize{$3.5318\textbf{E}-07$} \\
  \hline
  \hline
  \footnotesize{$1/4$}  & \footnotesize{$1/2000$}& \footnotesize{$2.7711\textbf{E}-05$} & \footnotesize{$1.8996\textbf{E}-05$} \\
  \footnotesize{$1/8$}  & \footnotesize{$1/2000$}& \footnotesize{$7.2794\textbf{E}-06$} & \footnotesize{$4.6964\textbf{E}-06$} \\
  \footnotesize{$1/16$} & \footnotesize{$1/2000$}& \footnotesize{$1.9737\textbf{E}-06$} & \footnotesize{$1.2544\textbf{E}-06$} \\
  \footnotesize{$1/32$} & \footnotesize{$1/2000$}& \footnotesize{$6.4648\textbf{E}-07$} & \footnotesize{$4.0951\textbf{E}-07$} \\
  \footnotesize{$1/64$} & \footnotesize{$1/2000$}& \footnotesize{$3.1454\textbf{E}-07$} & \footnotesize{$1.9998\textbf{E}-07$} \\
  \hline
  \hline
\end{tabular}
\end{center}
\end{table}
In addition, the results obtained for approximate solution along with estimate errors for time level $t=1, 5, 10, 20, 40$ and $t=100$ with $h=1/8$ and $t=1/1000$ in $-10 \leq x \leq 10$ are shown in Fig. 1.
The new method can be applied for Eq. (\ref{fn}) with other choices of $\rho$.
The space--time graph of approximate solution and related error estimate for $\rho=-1$ (the real Newell--Whitehead equation)
in domain $-10 \leq x \leq 10$ for times $t=0.25, 0.50, 0.75$ and $t=1.0$ are reported in Fig. 2.

\subsection{Example 2}

Consider the following generalized Fitzhugh–-Nagumo equation with time--dependent coefficients
\begin{equation} \label{10}
u_t+\cos(t)u_x-\cos(t)u_{xx}-2\cos(t)(u(1-u)(\rho-u))=0,\ \ \  (x,t) \in [a,b] \times [0,t].
\end{equation}
Also suppose that the initial and boundary conditions are taken from the exact solution given by \cite{Bhrawy,wazwazFN}
\begin{equation}\label{exact2}
u(x,t)=\frac{\rho}{2}+\frac{\rho}{2} \tanh(\frac{\rho}{2}(x-(3-\rho)\sin(t))).
\end{equation}
Table 2 lists the errors for this problem using the presented method with constant parameter $a=-b=1$, $\rho=1$, $\tau=0.001$, $t=1$ for different values of $h$.
\begin{table}[h]
\caption{The estimate errors for Ex. 2 for different values of $h$.}
\begin{center}
\begin{tabular}{ccc}
  \hline
  \hline
  h & $L_{\infty}-error$ & $RMS-error$ \\
  \hline
  \hline
  \footnotesize{$1/4$}   & \footnotesize{$1.0914\textbf{E}-03$} & \footnotesize{$9.5674\textbf{E}-04$} \\
  \footnotesize{$1/8$}   & \footnotesize{$3.4491\textbf{E}-04$} & \footnotesize{$2.9281\textbf{E}-04$} \\
  \footnotesize{$1/16$}  & \footnotesize{$1.5805\textbf{E}-04$} & \footnotesize{$1.2422\textbf{E}-04$} \\
  \footnotesize{$1/32$}  & \footnotesize{$1.1082\textbf{E}-04$} & \footnotesize{$8.2027\textbf{E}-05$} \\
  \footnotesize{$1/64$}  & \footnotesize{$9.8895\textbf{E}-05$} & \footnotesize{$7.1495\textbf{E}-05$} \\
  \footnotesize{$1/128$} & \footnotesize{$9.5897\textbf{E}-05$} & \footnotesize{$6.8814\textbf{E}-05$} \\
  \hline
  \hline
\end{tabular}
\end{center}
\end{table}
In addition, with the mentioned parameters, the obtained estimate errors for fixed $h=1/128$ and for various $\tau$ are reported in Table 3.
\begin{table}[h]
\caption{The estimate errors for Ex. 2 for different values of $\tau$}
\begin{center}
\begin{tabular}{ccc}
  \hline
  \hline
  $\tau$ & $L_{\infty}-error$ & $RMS-error$ \\
  \hline
  \hline
  \footnotesize{$1/100$}   & \footnotesize{$9.5923\textbf{E}-04$} & \footnotesize{$6.8834\textbf{E}-04$} \\
  \footnotesize{$1/200$}   & \footnotesize{$4.7752\textbf{E}-04$} & \footnotesize{$3.4244\textbf{E}-04$} \\
  \footnotesize{$1/400$}   & \footnotesize{$2.3862\textbf{E}-04$} & \footnotesize{$1.7109\textbf{E}-04$} \\
  \footnotesize{$1/800$}   & \footnotesize{$1.1965\textbf{E}-04$} & \footnotesize{$8.5833\textbf{E}-05$} \\
  \footnotesize{$1/1600$}  & \footnotesize{$6.0287\textbf{E}-05$} & \footnotesize{$4.3306\textbf{E}-05$} \\
  \footnotesize{$1/3200$}  & \footnotesize{$3.0634\textbf{E}-05$} & \footnotesize{$2.2068\textbf{E}-05$} \\
  \hline
  \hline
\end{tabular}
\end{center}
\end{table}
We also draw the approximate solutions along with estimate errors for six choices of $\rho=0.25, 0.50, 0.75, 1, 1.25$ and $\rho=1.5$ in $-10 \leq x \leq 10$ in Fig. 3. This figure shows that by increasing $\rho$ the absolute error increases. In addition, the approximate solution in conjunction with the exact solution for $\tau=0.001$, $\rho=1.5$ and $h=1/4$ at different times $t=0.25, 0.50, 0.75$ and $t=1$ are shown in Fig. 4.
\subsection{Example 3}

Consider the following generalized Fisher's equation
\begin{equation}\label{fisher}
u_t=u_{xx}+u(1-u^{\alpha}), \ \ \ (x,t) \in [-2,2] \times [0,1].
\end{equation}
The exact solution of the above equation is given by \cite{wazwazFisher,wangFisher}
\begin{equation}\label{fisheric}
u(x,t)=\{\frac{1}{2}\tanh[-\frac{\alpha}{2\sqrt{\alpha+4}}(x-\frac{\alpha+4}{\sqrt{\alpha+4}}t)+
\frac{1}{2}]+\frac{1}{2}\}^{2/\alpha}.
\end{equation}
We assume that the initial and the Dirichlet boundary conditions are taken from the above exact solution.
The approximate solutions and the estimate errors for this problem when $\alpha=1, 2, 3, 4, 5$ and $\alpha=6$ are shown in Fig. 5.
As this figure shows, the estimate errors increase when the value of $\alpha$ increases.

\section{Conclusion}

In this article, a numerical method based on the dual reciprocity boundary elements method (DRBEM) is outlined for solving the one--dimensional nonlinear parabolic partial differential equations. The list of equations investigated includes Fisher's equation, generalized Fisher's equation, Allen-–Cahn equation, Newell-–Whithead equation, Fitzhugh–-Nagumo equation and  generalized Fitzhugh-–Nagumo equation with time variable coefficient.
The dual reciprocity idea was applied to eliminate the domain integrals appearing in the boundary integral equation.
Linear radial basis functions (RBFs) were used in the presented method as approximate functions.
We used the implicit finite difference method in time and the boundary integral equation technique in space to
discretize the main differential equation and convert it to a linear algebraic system of equations. The
nonlinear terms are treated iteratively within each time step by using a simple predictor--corrector scheme. Numerical
results are presented for some test problems to demonstrate the usefulness and accuracy of the
new method.

\section{conflict of interest statement:}

On behalf of all authors, the corresponding author states that there is no conflict of interest.

\section{Data Availability Statement:}

The data used in this research article is not applicable as no specific data sets were utilized. The conclusions and findings presented in this paper are based on theoretical analysis, literature review, and other relevant scholarly resources. All references cited are available in the reference section for further examination.

\begin{small}

\end{small}

\newpage
\begin{figure}[t]
\includegraphics[width=8cm, height=7.cm]{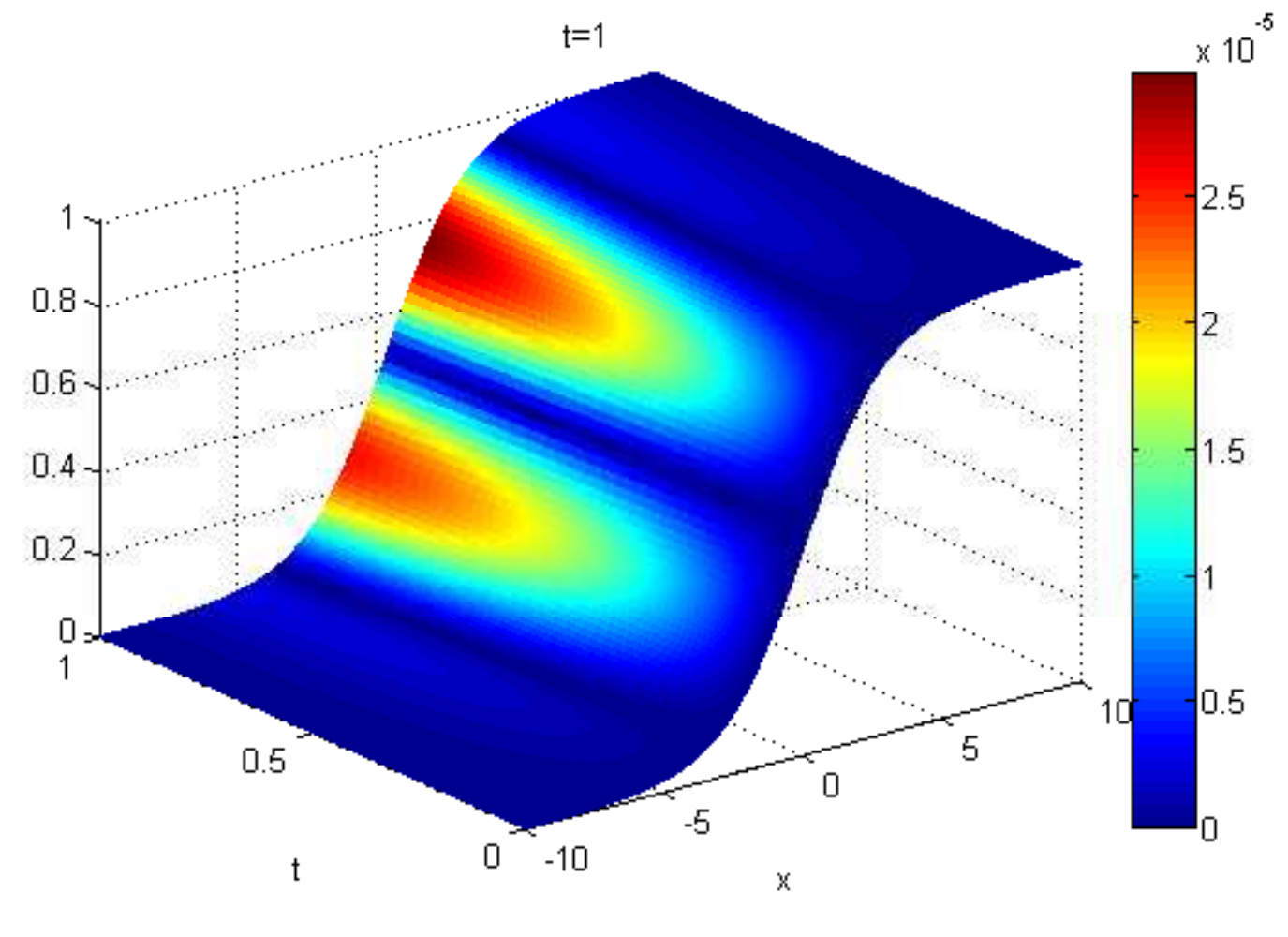}
\includegraphics[width=8cm, height=7.cm]{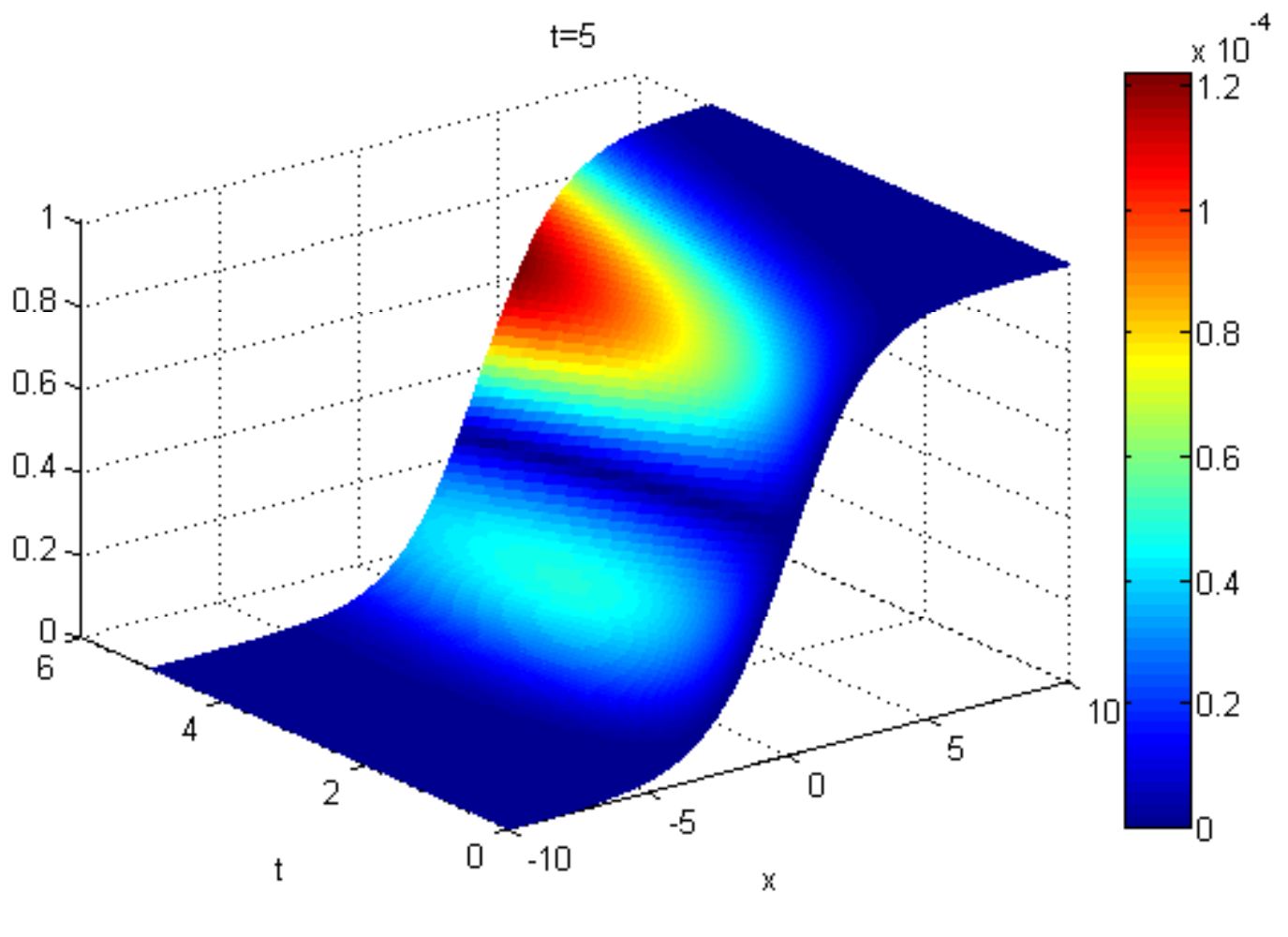}
\includegraphics[width=8cm, height=7.cm]{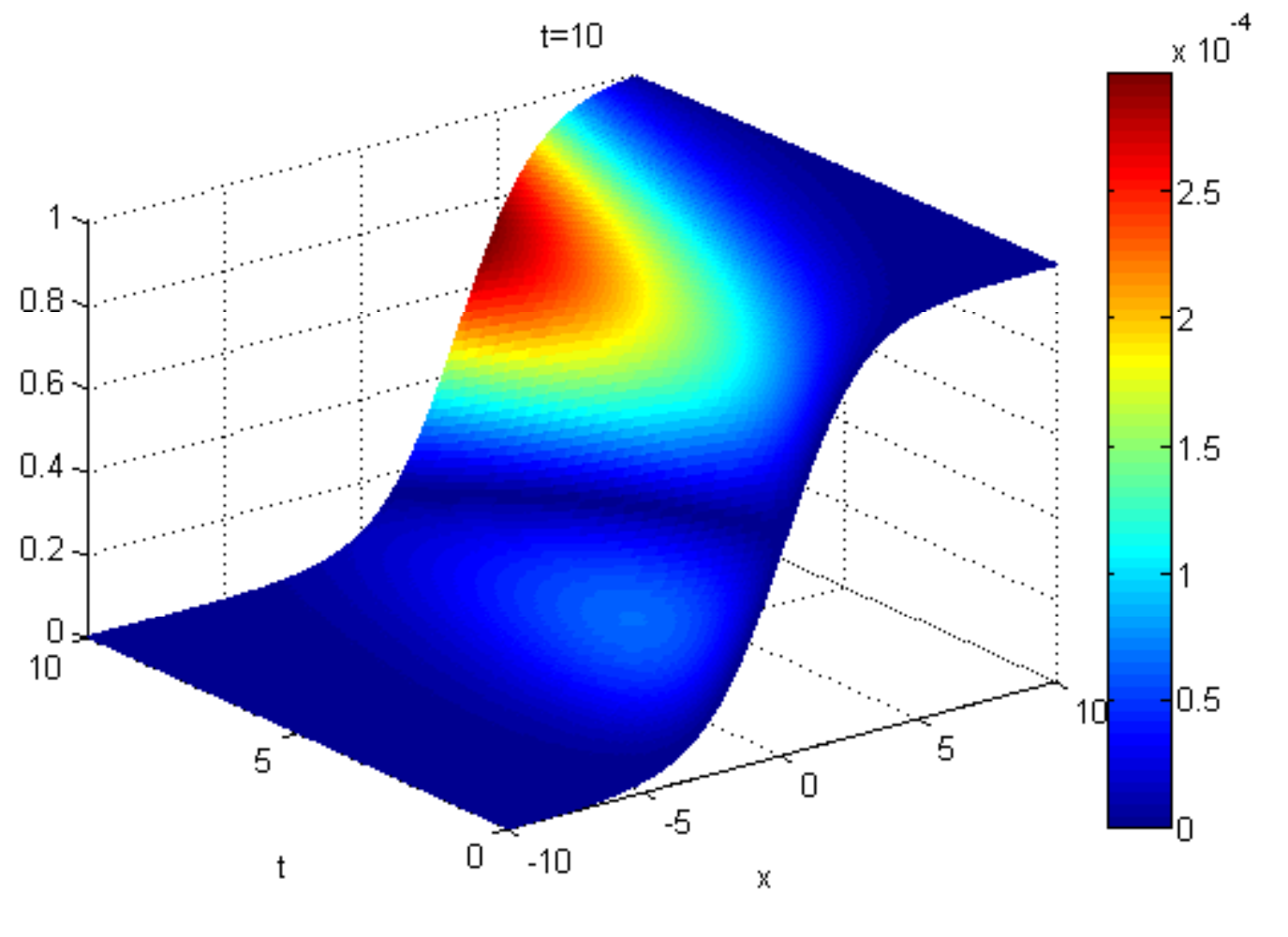}
\includegraphics[width=8cm, height=7.cm]{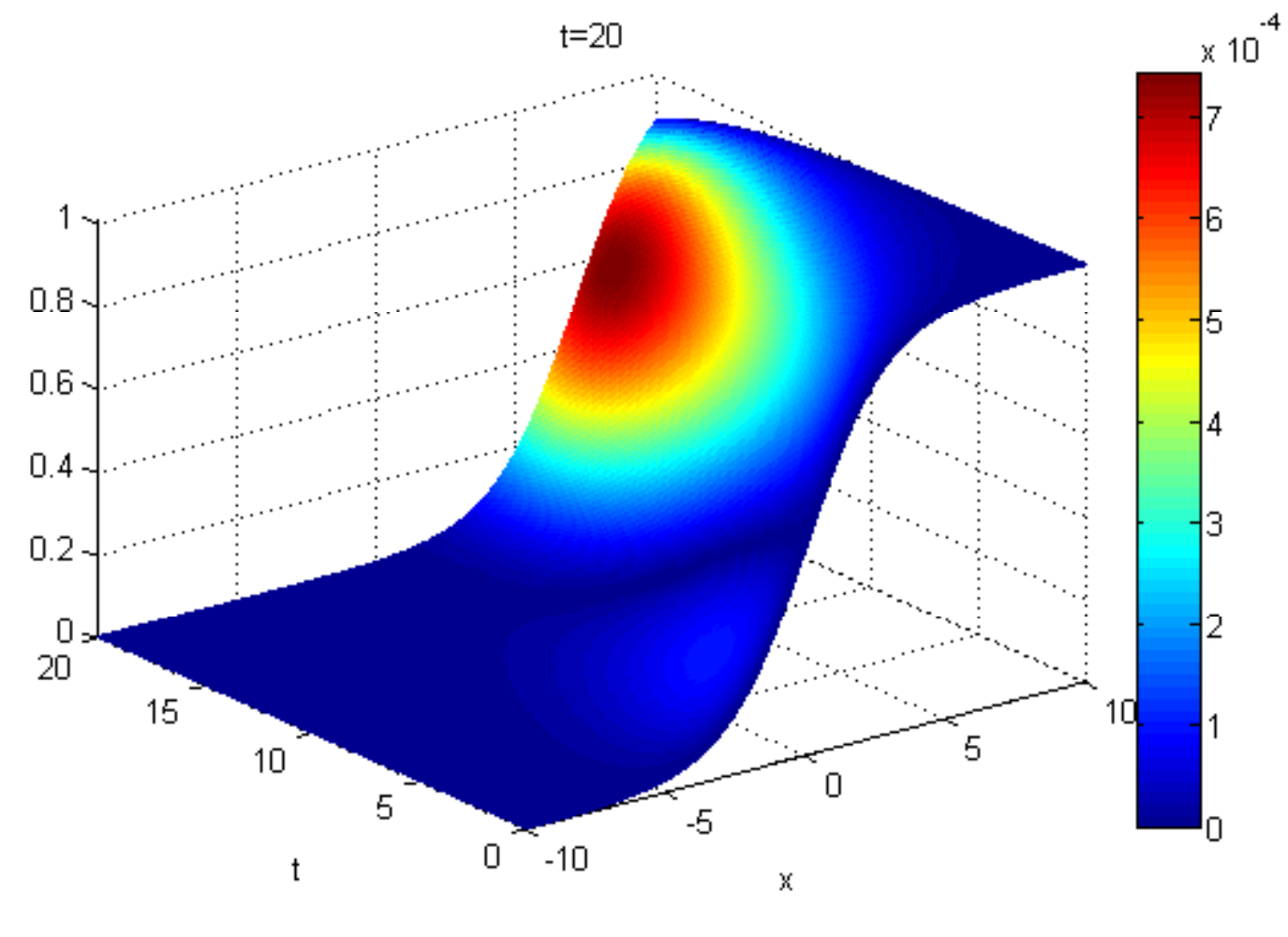}
\includegraphics[width=8cm, height=7.cm]{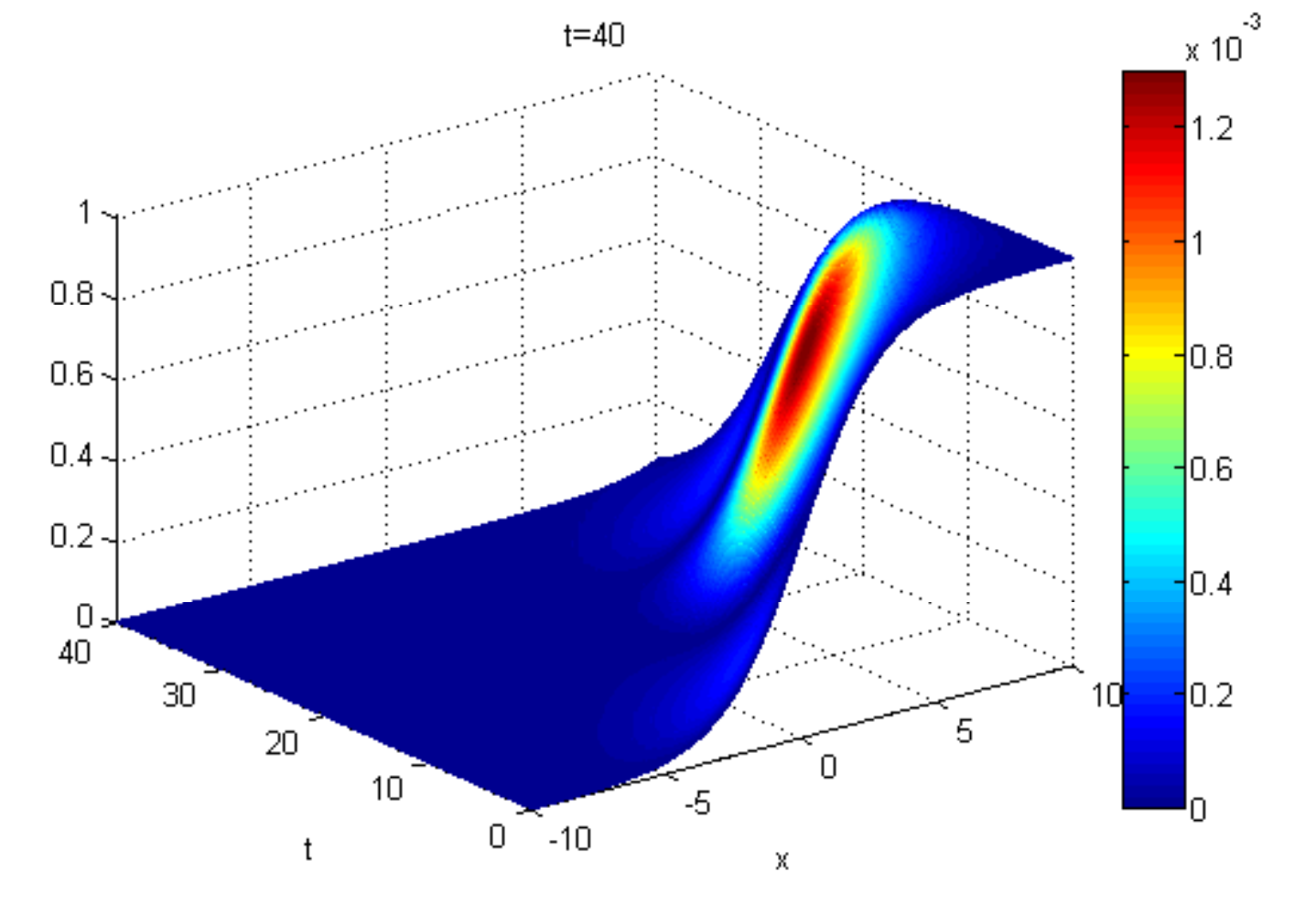}
\includegraphics[width=8cm, height=7.cm]{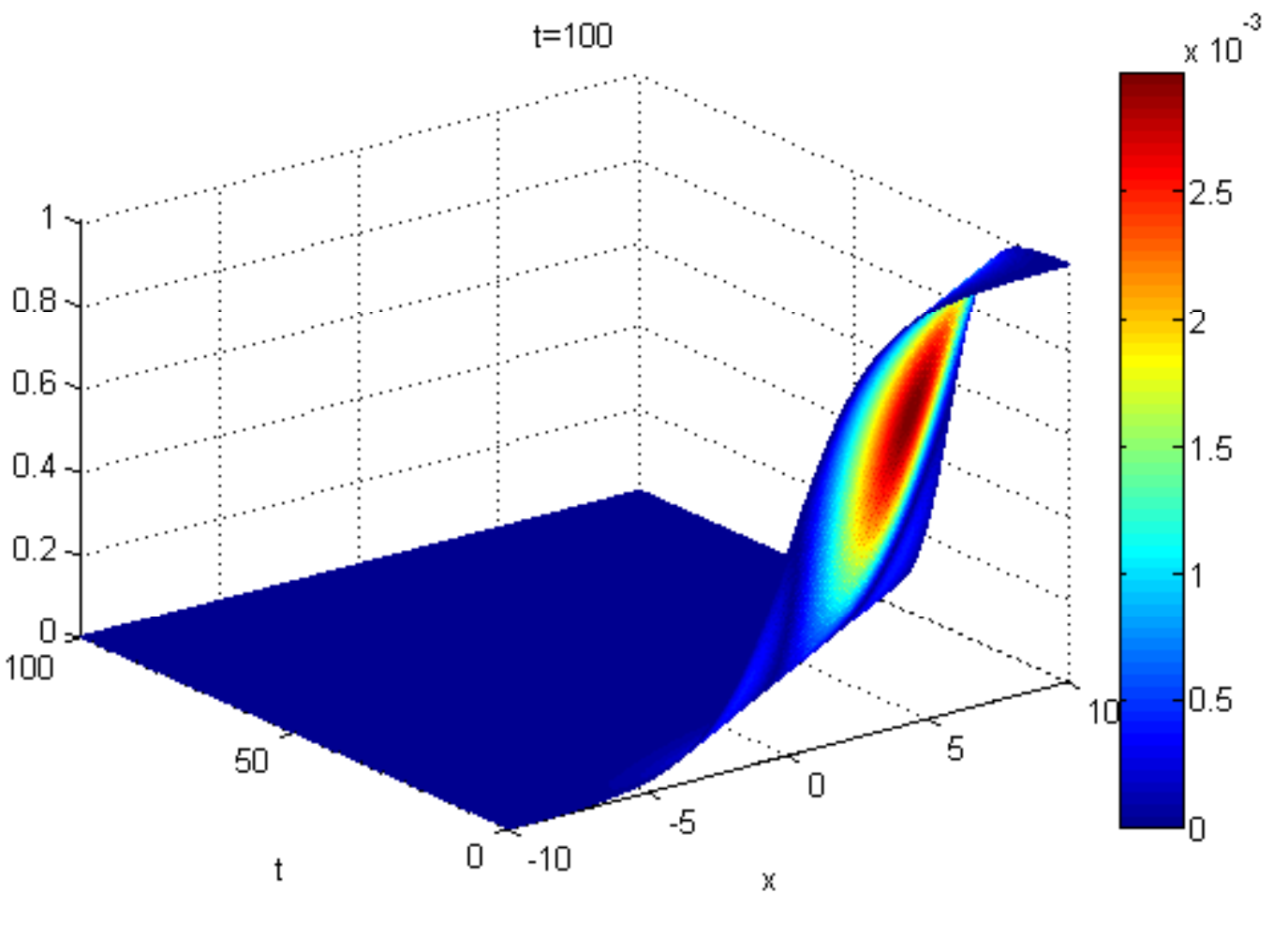}
\caption{Graphs of approximated solutions along with estimate errors at times $t=1, 5, 10, 20, 40$ and $t=100$ obtained for
Example 1 for constants $t=1/1000$, $h=1/8$ and $\rho=3/4$ in $-10 \leq x \leq 10$.}
\end{figure}

\begin{figure}[t]
\includegraphics[width=8cm, height=7.cm]{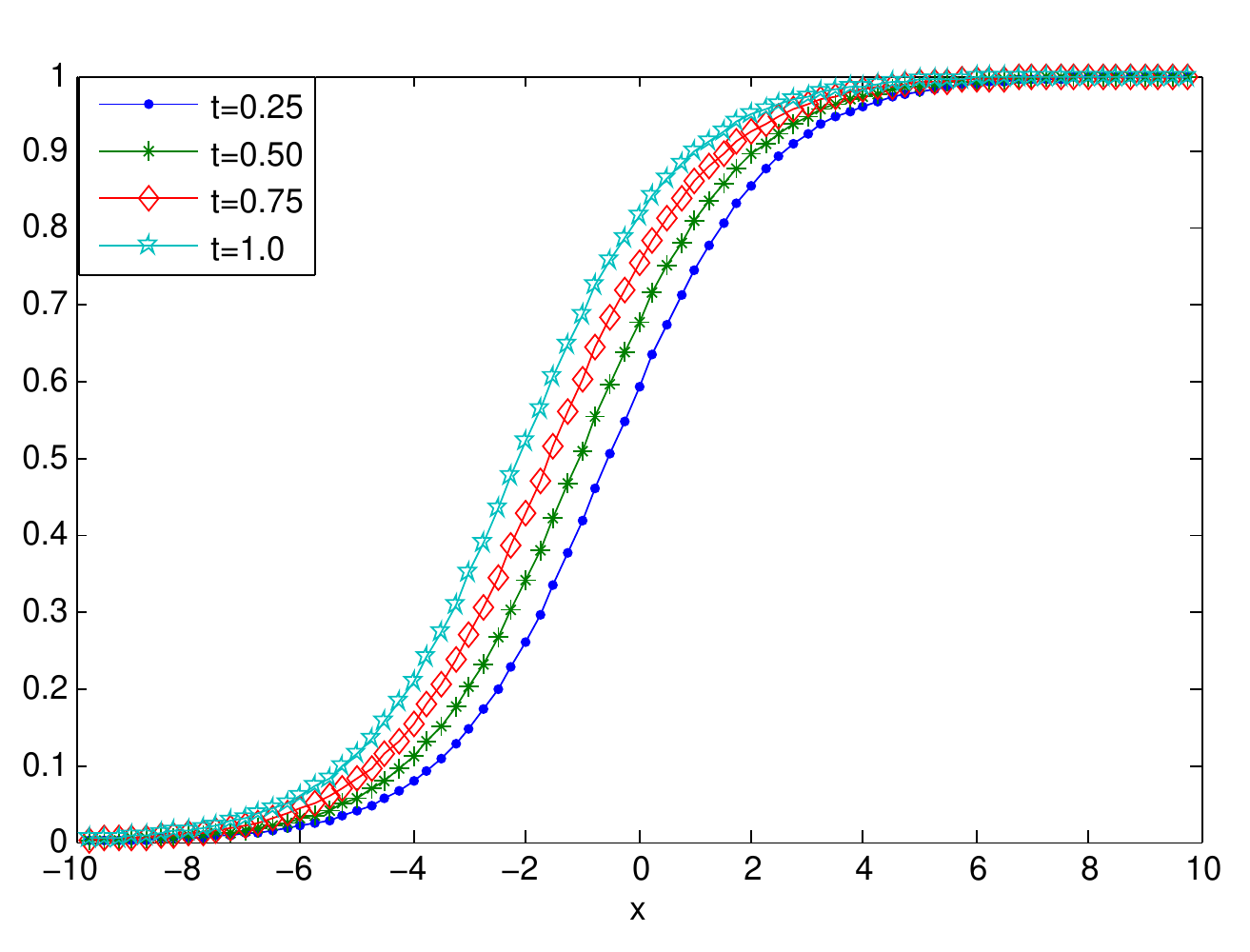}
\includegraphics[width=8cm, height=7.cm]{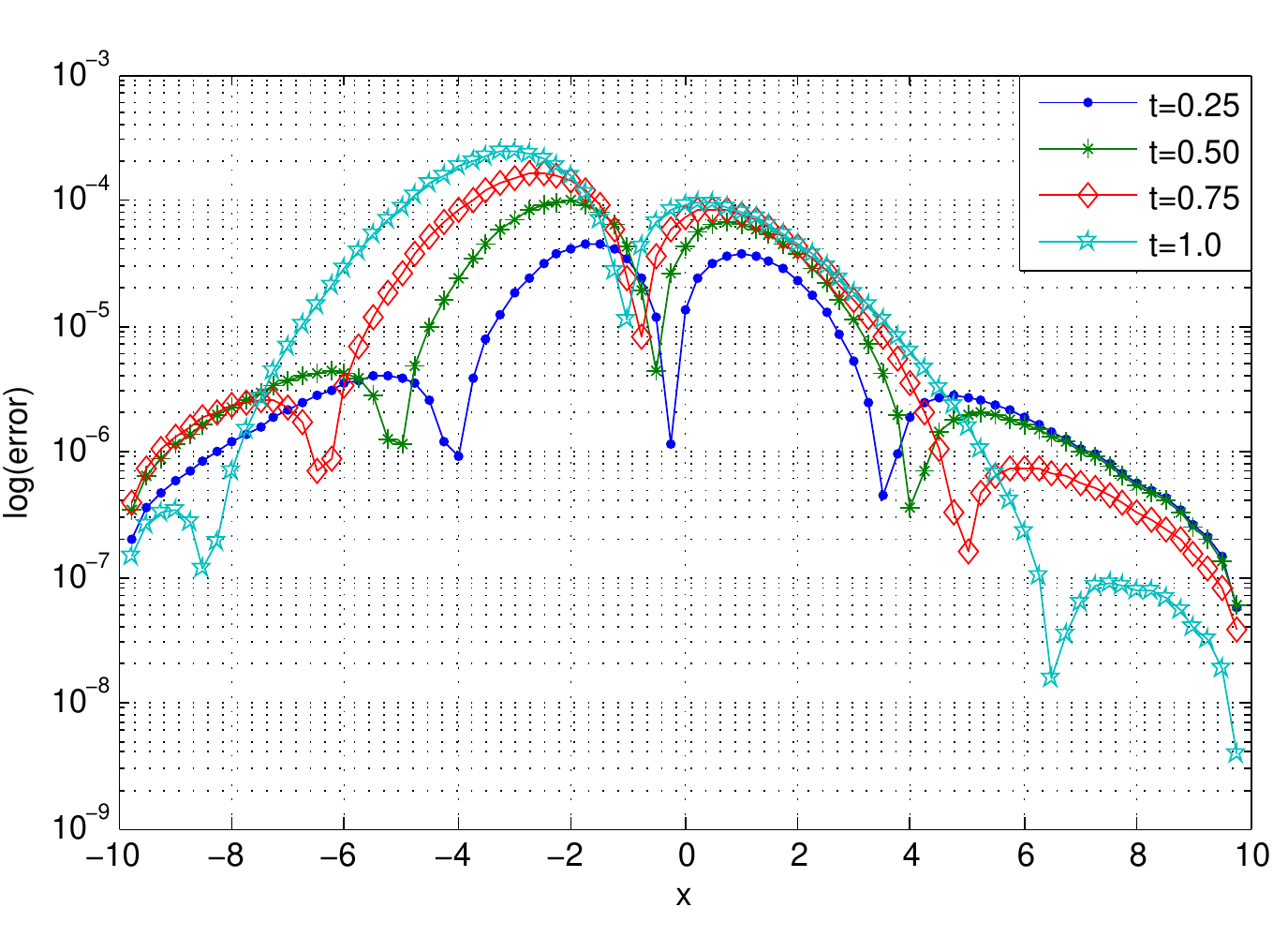}
\caption{Graphs of approximated solutions (left plan) and their estimate errors (right) at times $t=0.25, 0.5, 0.75$ and
 $t=1.0$ obtained for Example 1 for constants $\tau=1/1000$, $h=1/8$ and $\rho=-1$ in $-10 \leq x \leq 10$.}
\end{figure}


\begin{figure}[t]
\includegraphics[width=8cm, height=7.cm]{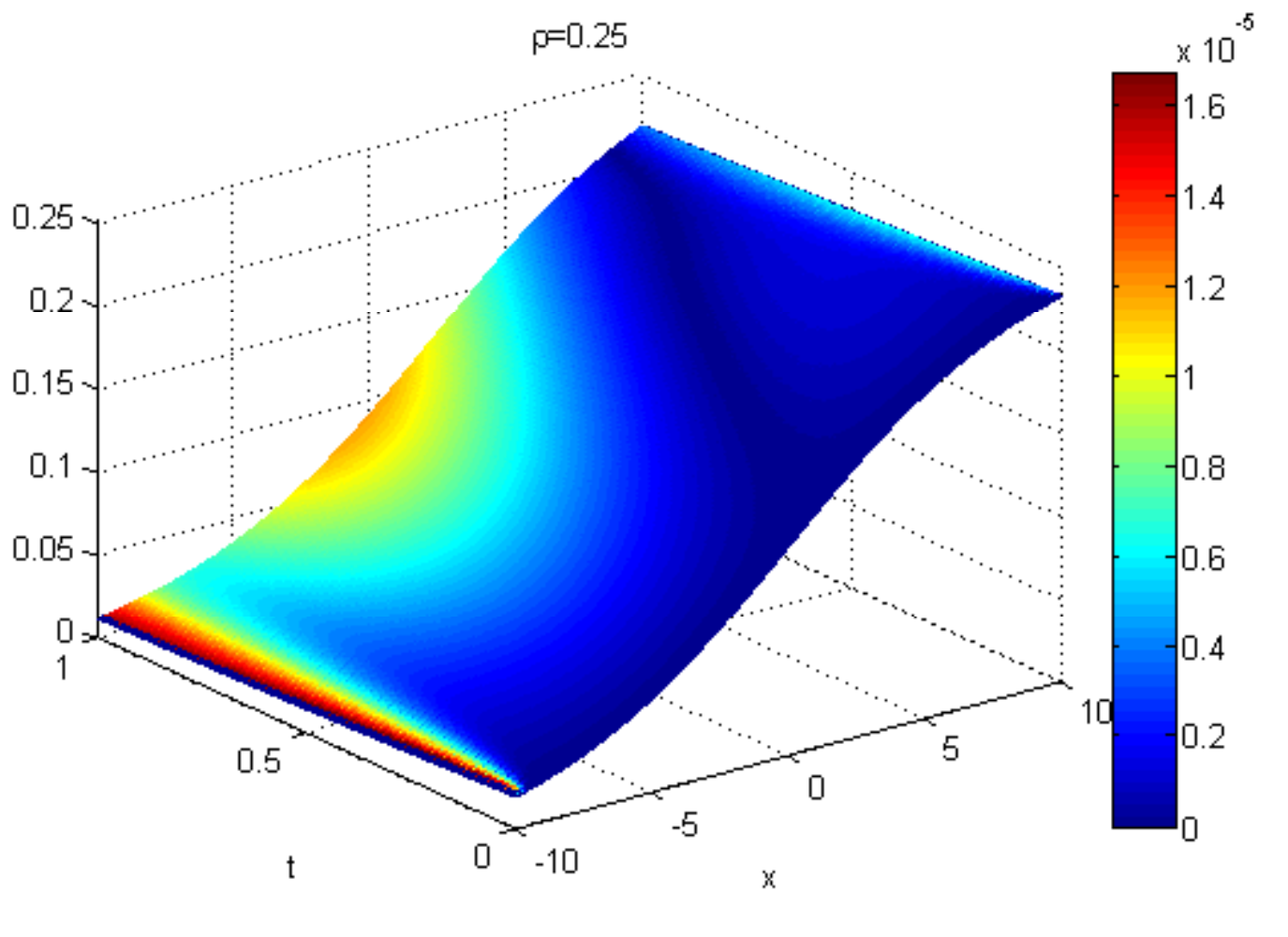}
\includegraphics[width=8cm, height=7.cm]{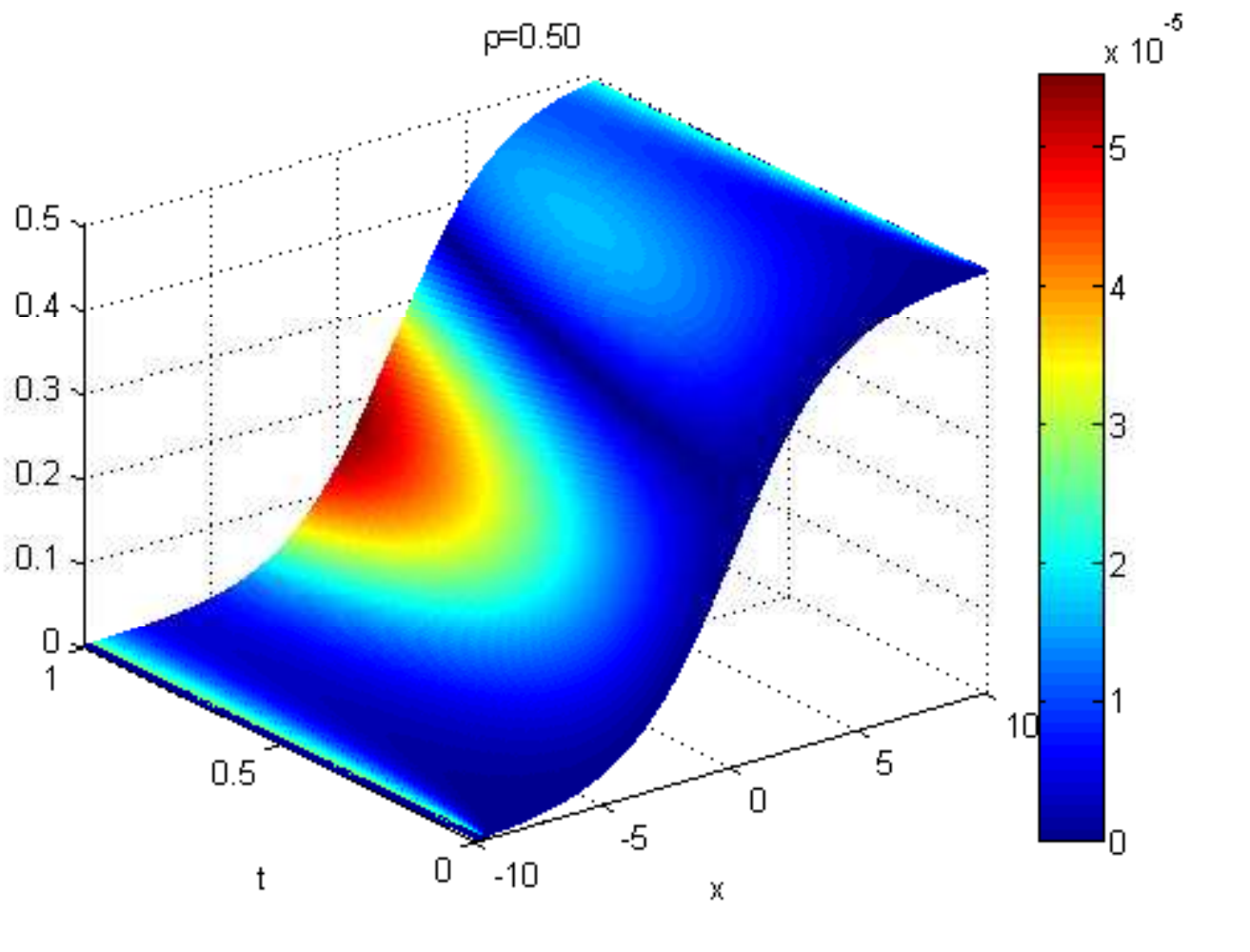}
\includegraphics[width=8cm, height=7.cm]{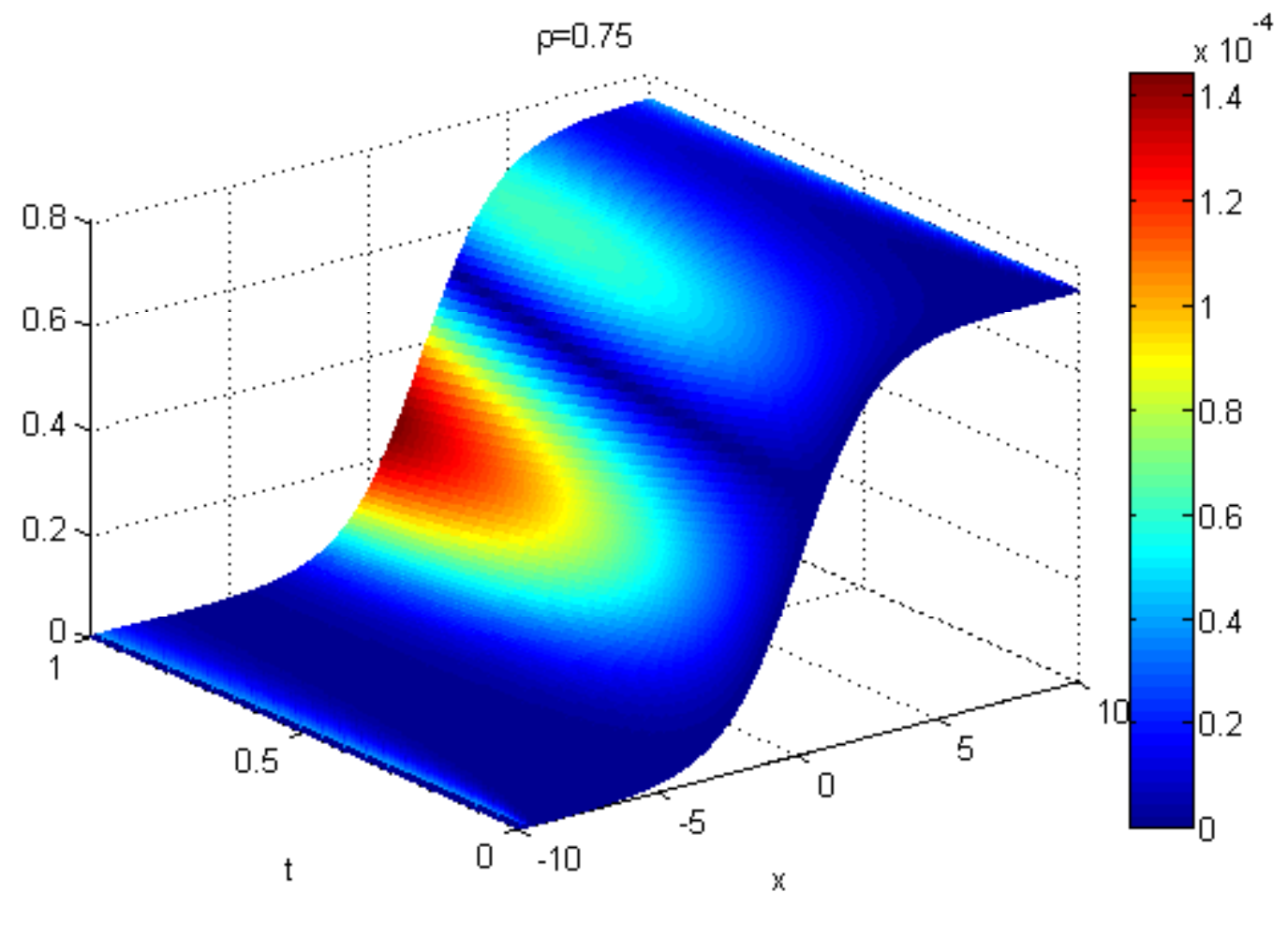}
\includegraphics[width=8cm, height=7.cm]{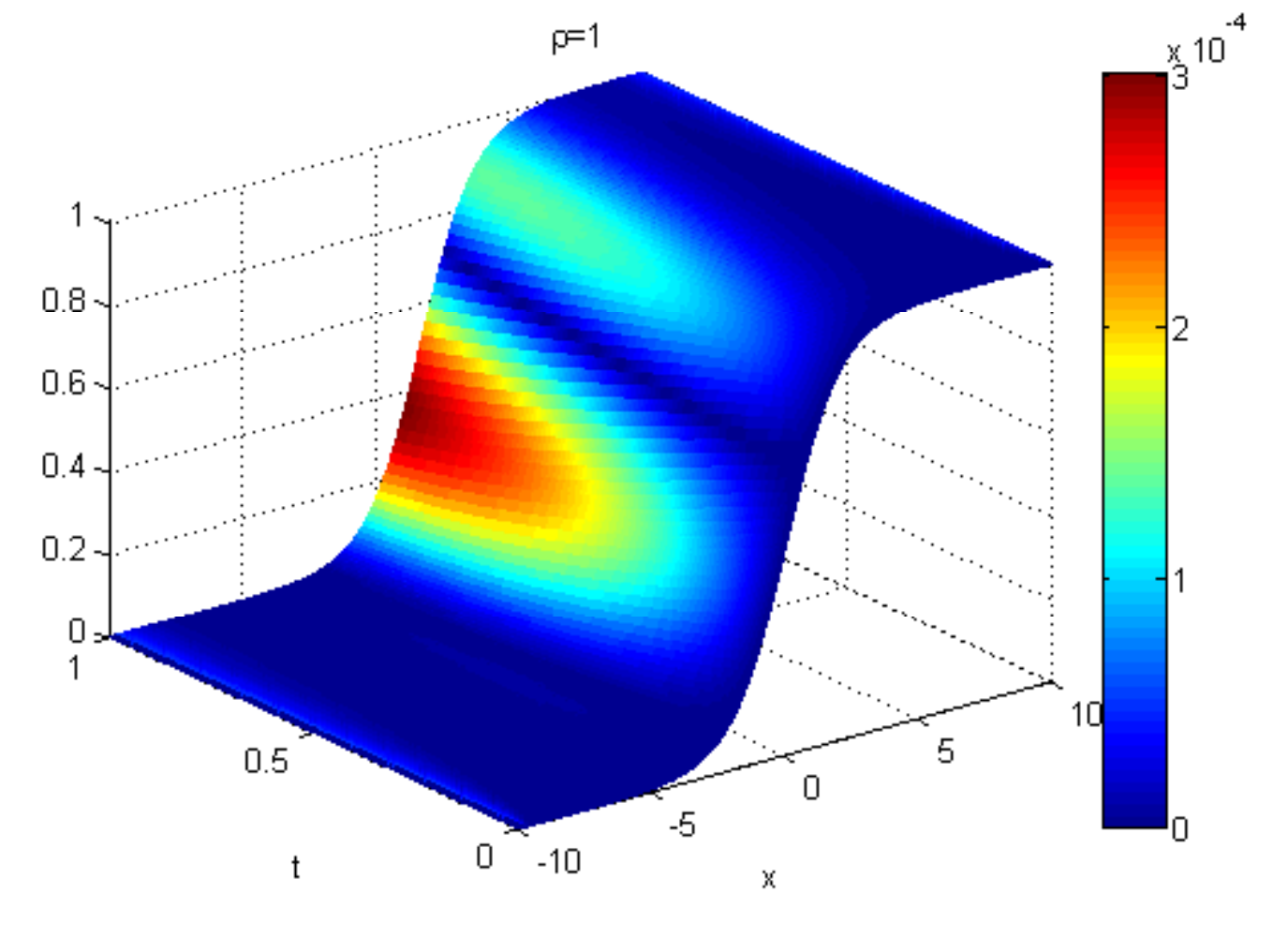}
\includegraphics[width=8cm, height=7.cm]{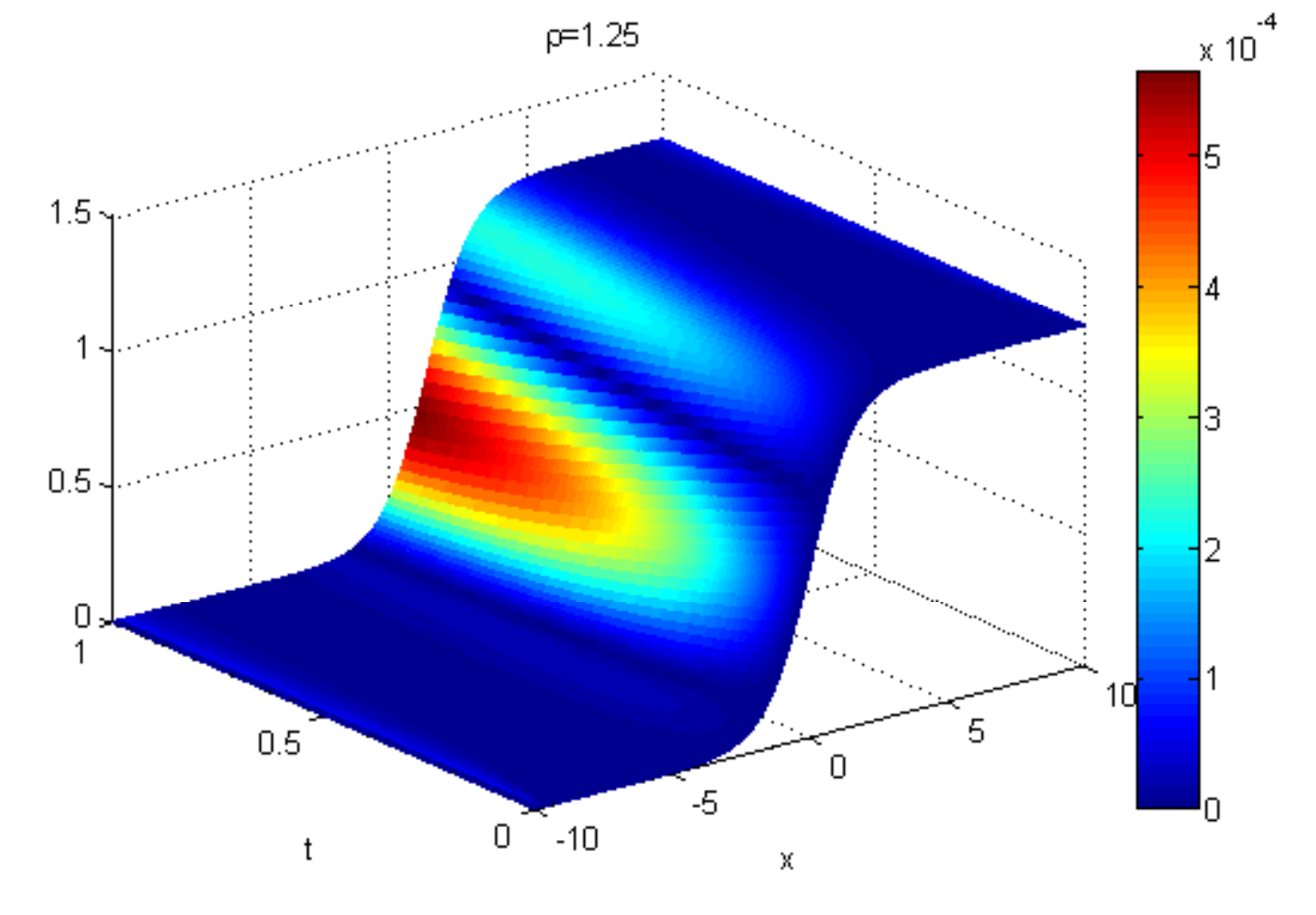}
\includegraphics[width=8cm, height=7.cm]{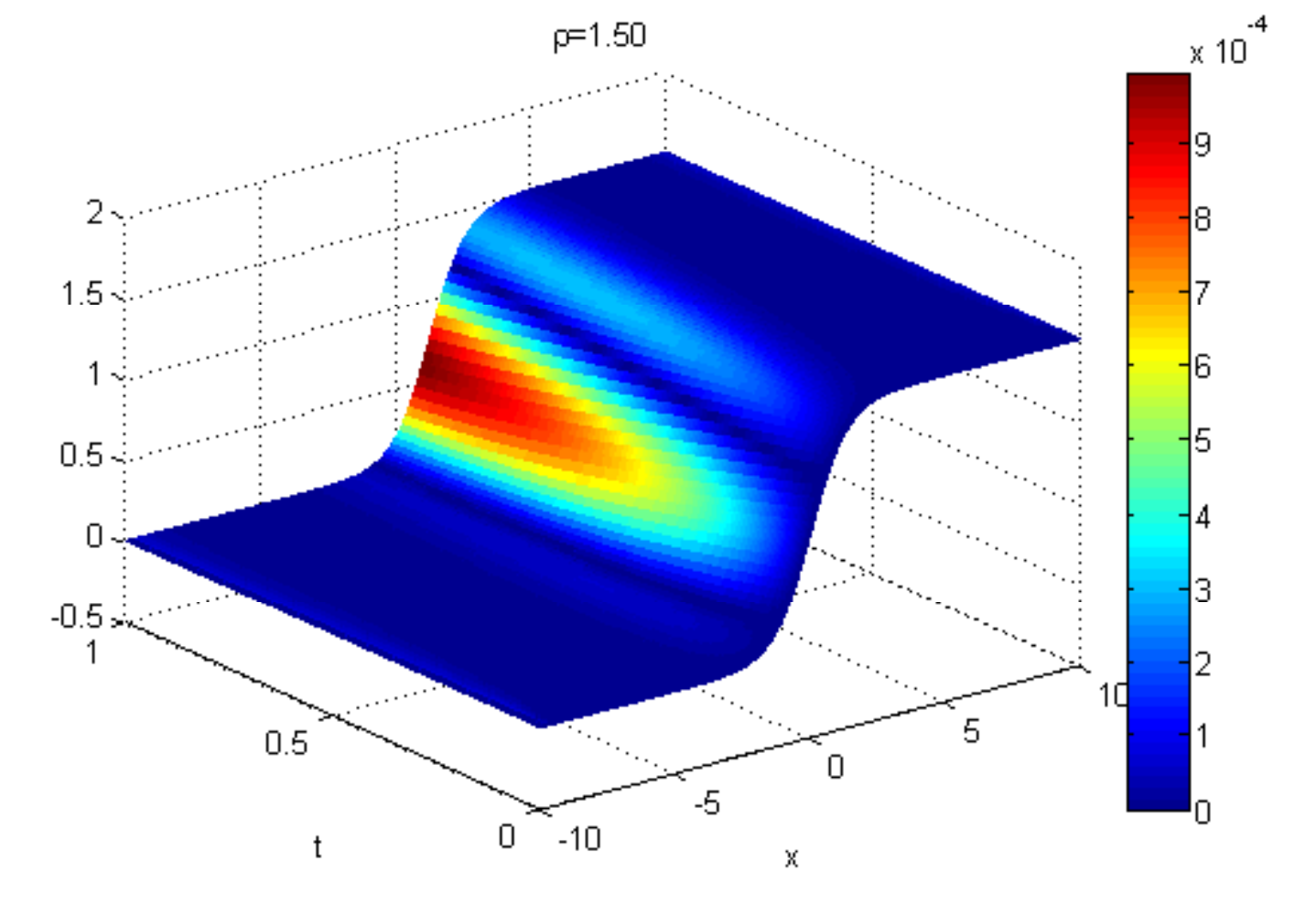}
\caption{Graphs of approximated solutions along with estimate errors at times $t=1$ for various $\rho=0.25, 0.50, 0.75, 1, 1.25$
and $\rho=1.50$ obtained for Example 2 for constants $t=1/1000$ and $h=1/8$ in $-10 \leq x \leq 10$.}
\end{figure}

\begin{figure}[t]
\begin{center}
\includegraphics[width=8cm, height=7.cm]{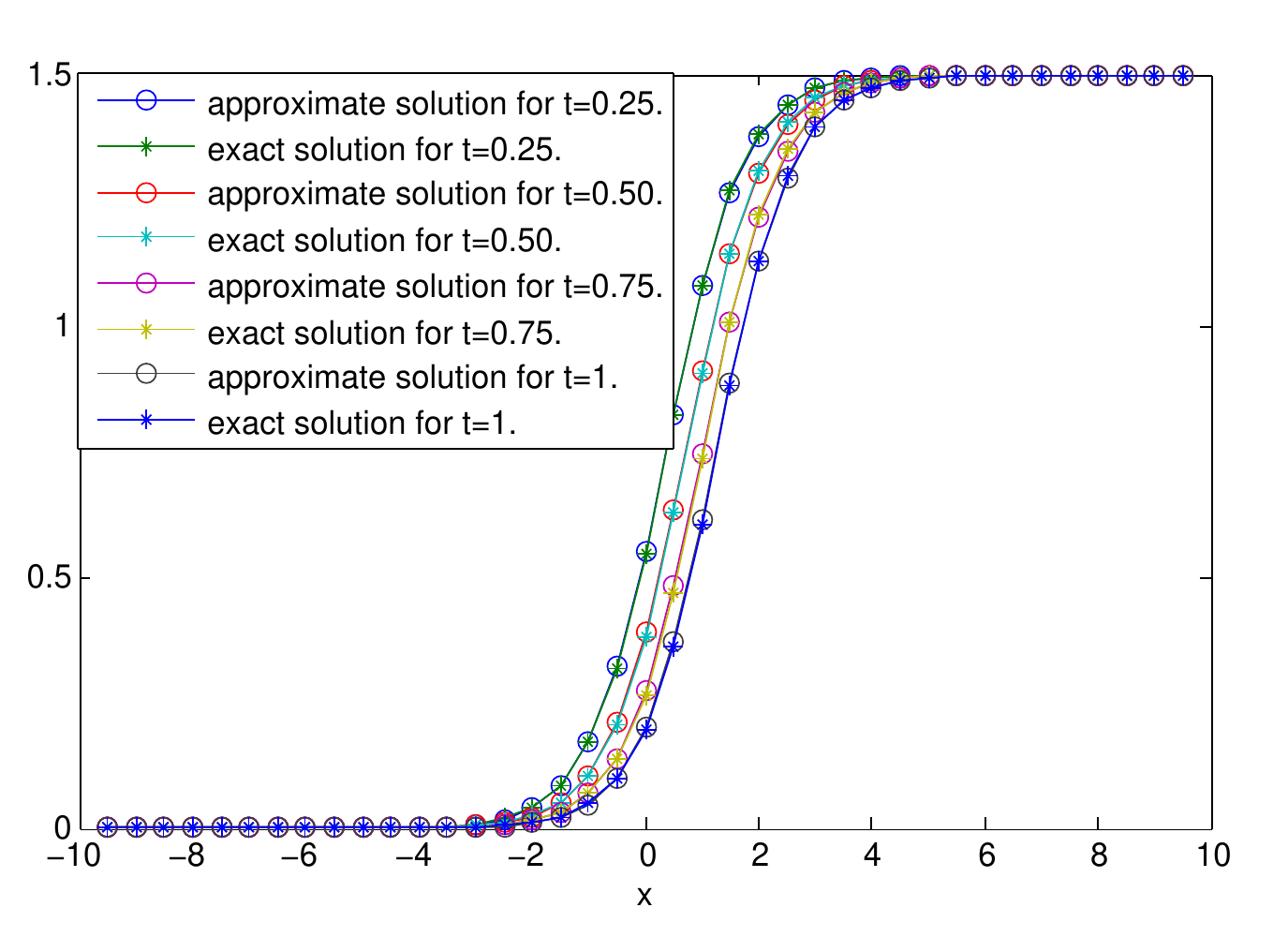}
\caption{Graphs of approximated solutions along with exact solutions at times $t=0.25, 0.50, 0.75$ and $t=1$
with $\rho=1.50$  for Example 2 for constants $t=1/1000$ and $h=1/8$ in $-10 \leq x \leq 10$.}
\end{center}
\end{figure}


\begin{figure}[t]
\includegraphics[width=8cm, height=7.cm]{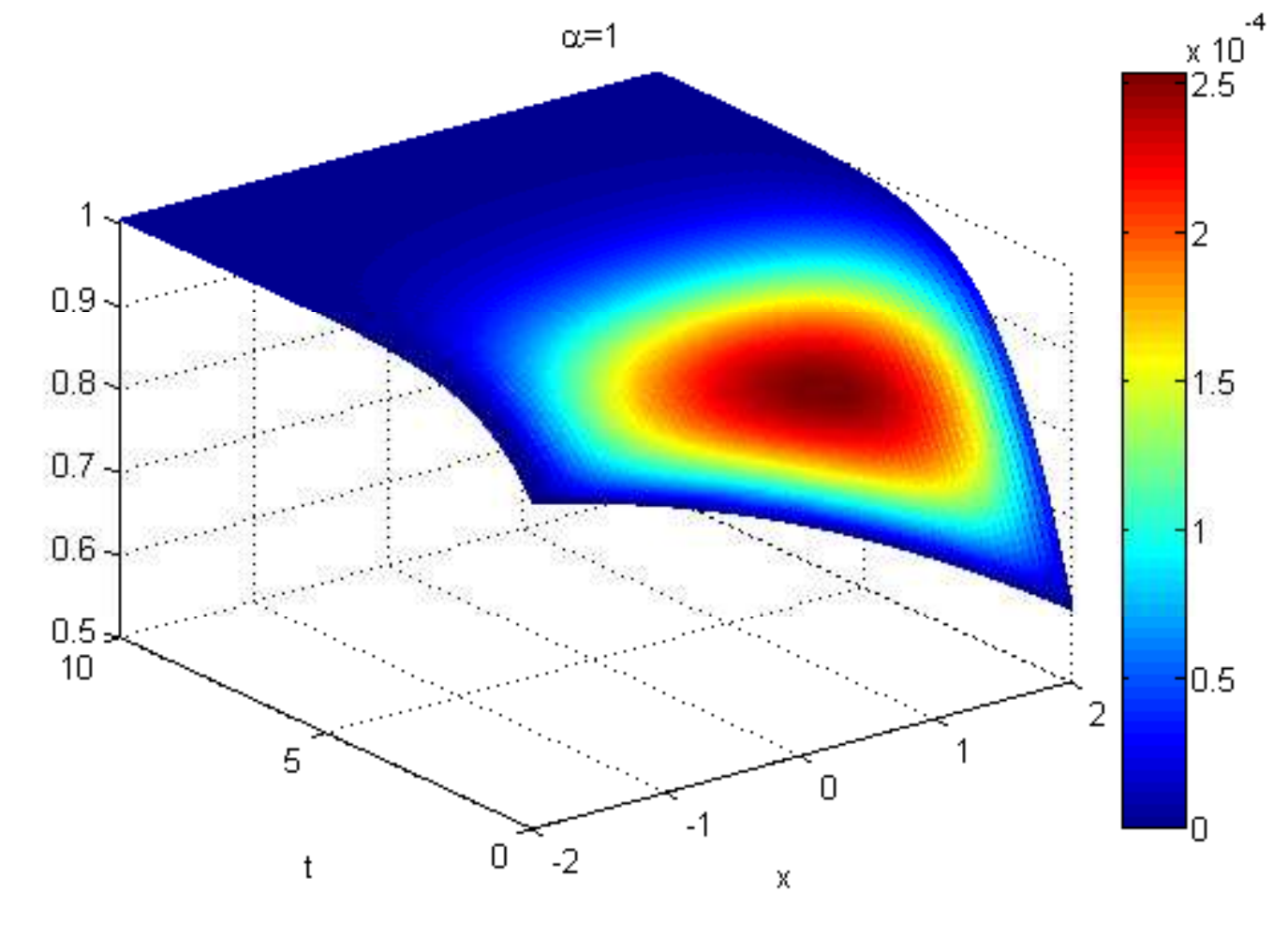}
\includegraphics[width=8cm, height=7.cm]{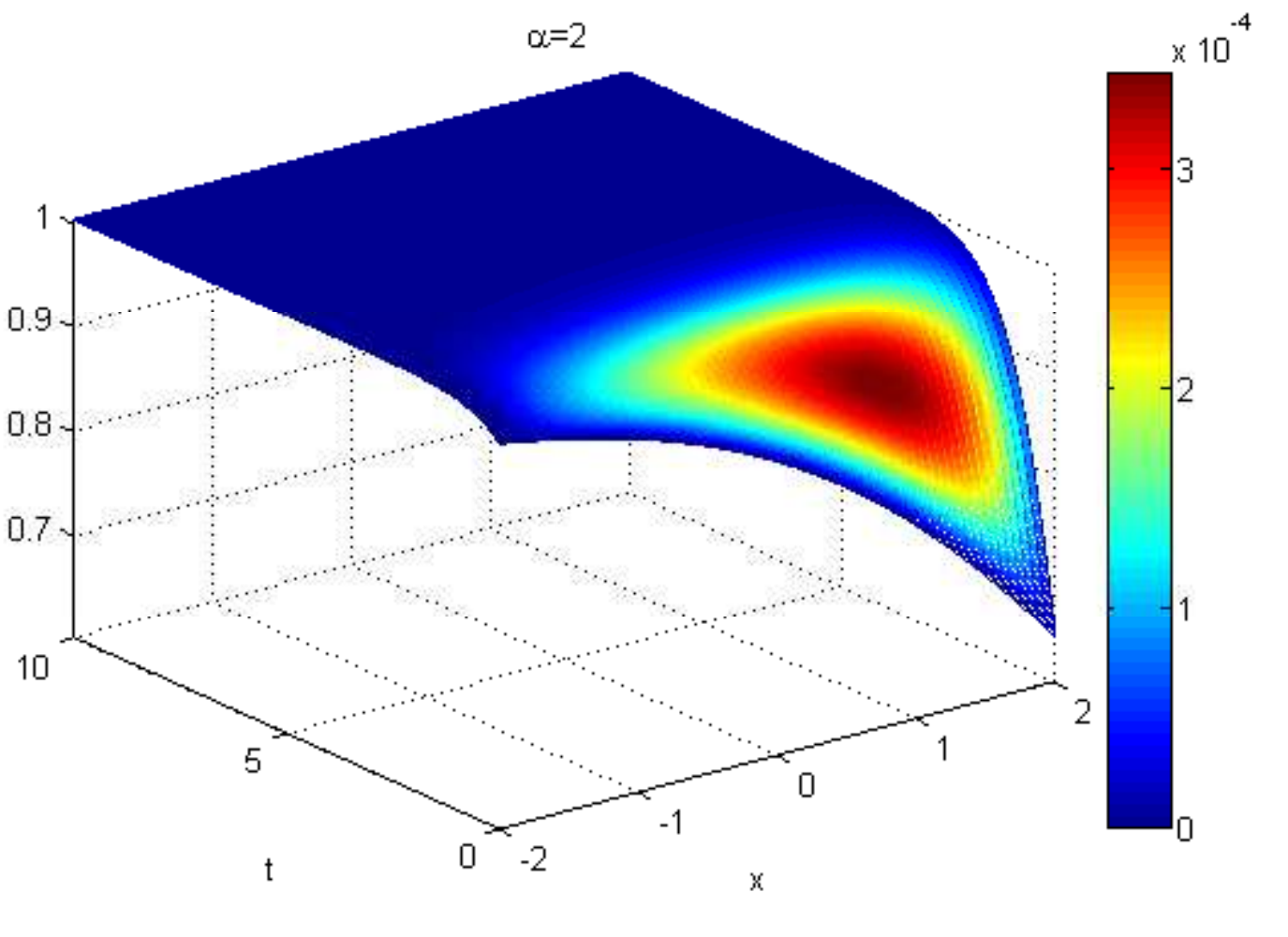}
\includegraphics[width=8cm, height=7.cm]{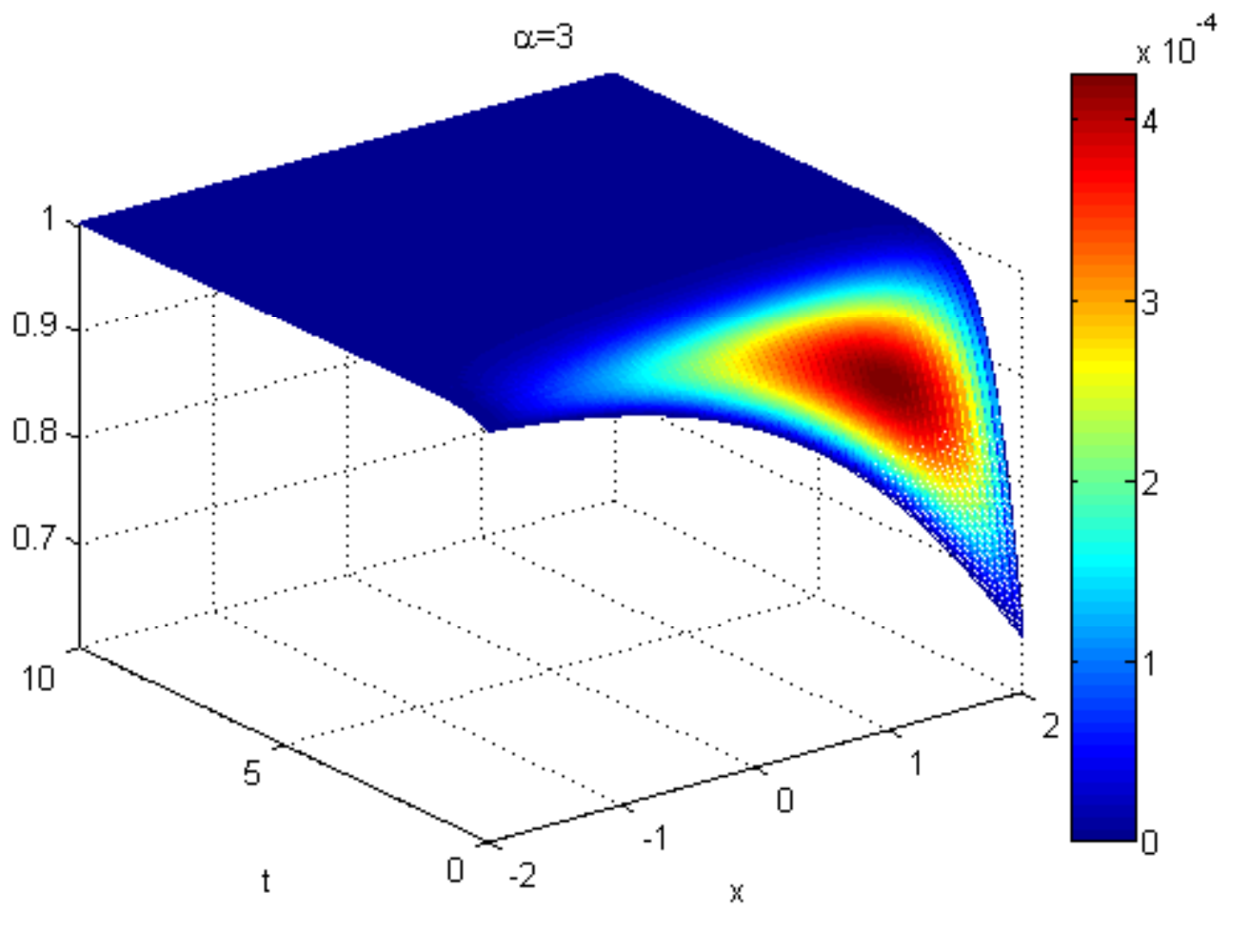}
\includegraphics[width=8cm, height=7.cm]{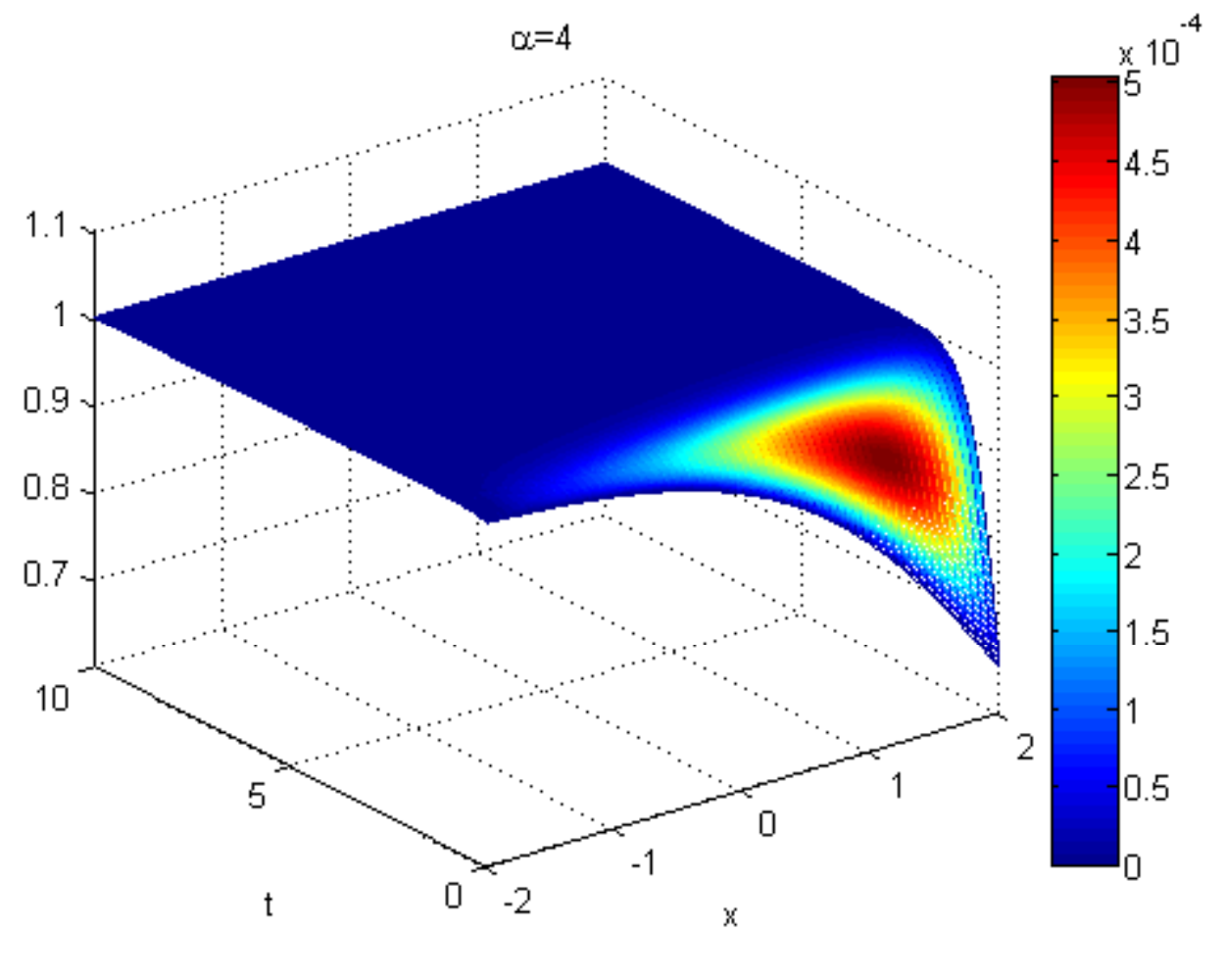}
\includegraphics[width=8cm, height=7.cm]{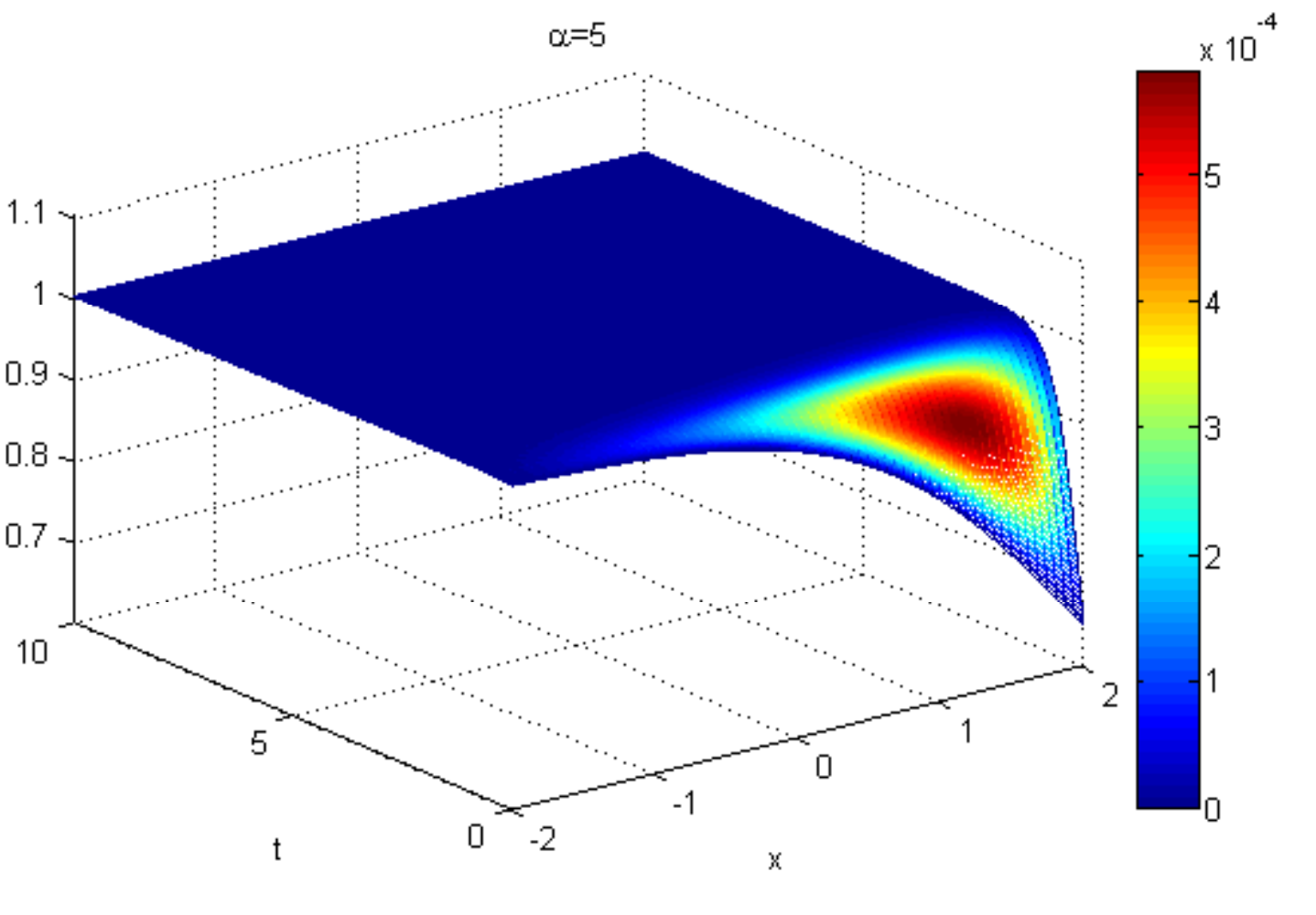}
\includegraphics[width=8cm, height=7.cm]{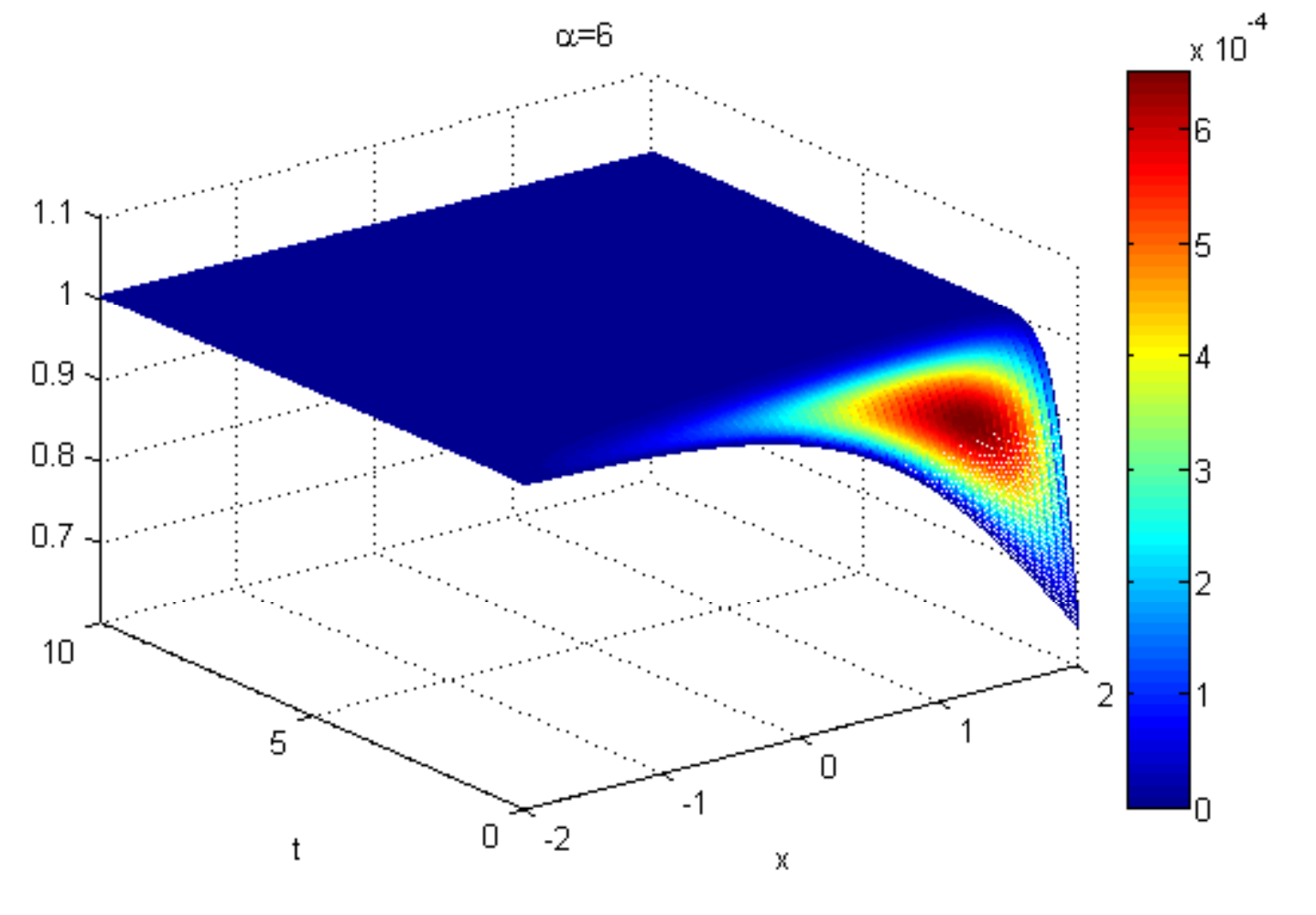}
\caption{Graphs of approximated solutions along with estimate errors for various $\alpha=1, 2, 3, 4, 5$ and $\alpha=6$
at time $t=10$ obtained for generalized Fisher's equation for constants $t=1/1000$ and $h=1/16$ in $-2 \leq x \leq 2$.}
\end{figure}

\end{document}